%% file: filtered_hyperinterp_shell_arxiv_ver2.tex
\documentclass[a4paper,11pt]{article}
\input{my_arxiv_header.tex}
\title{}
\author{Yoshihito Kazashi}
\bibliography{filtered_abbrev}

\newcommand{\sumQN}{\ensuremath{\sum_{j=0}^{\nu_Q}}}

\title{A fully discretised filtered polynomial approximation on spherical shells 
}
\begin{document}
\maketitle
\begin{abstract}
 A fully implementable filtered polynomial approximation on spherical shells is considered. The method proposed is a quadrature-based version of a filtered polynomial approximation. The radial direction and the angular direction of the shells are treated separately with constructive filtered polynomial approximation. The approximation error with respect to the supremum norm is shown to decay algebraically for functions in suitable differentiability classes. Numerical experiments support the results.
\end{abstract}

\keywords{Filtered polynomial approximation, Filtered hyperinterpolation, Spherical shells}
%\MSC[2010] 41A10 \sep 65D10

%\end{frontmatter}

%\linenumbers

\section{Introduction}
This paper is concerned with constructive global polynomial approximation on spherical shells. 
Problems on such domains naturally arise in a wide range of geosciences, and numerous computational methods are proposed
\cite{Flyer.N_etal_2014_shell,Fornberg.B_Flyer_2015_Acta,McWilliams.J.C_1996_shell,Sakuraba.A_dynamo,Tackley.P_2008_mantle_convection}. Nonetheless, theoretical analysis does not seem to have attracted much attention. One recent result considered in \cite{Kazashi.Y_2016_nonfiltered} is a fully discretised polynomial approximation on spherical shells. The method considered there can be seen as an approximation of the $L^2$-orthogonal projection. However, $L^2$-projection is not the best choice when one wants a small point-wise error---recall how the Fourier series of $f$ on the torus may fail to converge on any measure zero set if $f$ is merely continuous \cite{Kahane.J_Katznelson_1966_Fourier_divergence}.
 This paper considers a method with good uniform convergence using {\it filtering}. Our ultimate goal is to construct a fully discretised filtered polynomial approximation method and analyse the errors.

A classical remedy for the failure of the Fourier series on the torus mentioned above is to use a smoothing (or {\it filtering}) process such as Ces\'{a}ro sums, Lanczos smoothing, or the raised cosine smoothing \cite{Canuto.C2006book}. 
An underlying idea is to smoothly truncate the series by multiplying the Fourier coefficients of higher order by a suitably small factor. 

\red{Filtered approximations have also been considered for other regions, including the sphere and the ball. A standard approach to show the convergence rate is to show the uniform boundedness of the family of linear operators defined by the approximation, which readily reduces the problem to the best polynomial approximation rate. To show such boundedness, the assumption that the Fourier coefficients are given exactly is often fully utilised. See \cite{Petrushev.P_Xu_2005_interval_Jacobi,Sloan.I.H_2011_dela_Vallee,Dai.F_2011_Sobolev,Dai.F2013book} and references therein.}

In most realistic applications the Fourier coefficients are not known, as integrals are not computable exactly. As an alternative, quadrature-based approximations of these filtered methods have been considered for various settings, particularly on the sphere \cite{Mhaskar.H.N_Prestin_2005_fast_tour,Sloan.I_Womersley_2012_filtered_hyperinterp}. The present paper considers a quadrature-based filtered polynomial approximation for spherical shells $\shell=\{\bs{x}\in\mathbb{R}^3\ |\ \rin\le \norm{\bs{x}}_2\le\rout=\rin+\ep\}$, where $\rin\le 1\le \rout$ and $\rout-\rin=\ep>0$ as a domain. 
That is, given a function $f$, we smoothly truncate its Fourier series by a suitable `filter' $h$, and approximate Fourier coefficients by quadrature rules. The motivation is to propose an implementable technique with a good point-wise convergence. Our results give, to the best of our knowledge, the first theoretical results on constructive filtered polynomial approximation on spherical shells. Our method requires only point values of the function we approximate, and thus can be implemented exactly in a real number model of computation.

We regard $\shell$ as a product of the interval $\rinout$ in the radial direction and the unit sphere $S^2$ in the angular direction. % We consider approximations separately on $\rinout$ and $S^2$, in particular, we propose a constructive polynomial approximation.
 The product setting is natural, since in practice functions on a spherical shell vary on different scales in the radial and angular directions. For example, the mantle can be seen as a set of spherical layers with different characteristics \cite{Bangerth.W_2016_mantle_SIAM_news,Jackson_book_mantle,Davies_book_dynamic_earth}. Some properties of the atmosphere, such as the ionisation rate \cite[p. 151]{Earth_electrical_environment}, electric field \cite[p. 155]{Earth_electrical_environment}, depend strongly on the altitude, and hence vary rapidly in the radial direction.

%The $L^2$-orthogonal projection onto a space of orthogonal polynomials is a natural way to approximate functions. However, typically Fourier coefficients can be obtained only approximately, due to the necessity of evaluating the integrals. Instead, we employ methods in the radial and the angular direction that shares the same properties with $L^2$-projections. More specifically, f
For a continuous function $f\in C(\shell)$ we consider the approximation taking the following form. Let $J_k^{(\alpha,\beta)}$ $(\alpha,\beta>-1)$ be the Jacobi polynomial of degree $k$ mapped affinely to $\rinout$ from $[-1,1]$, and $Y_{\ell m}$ be the spherical harmonics of degree $\ell$ and order $m$. More detailed definitions are given later. Let $h\from[0,\infty)\times[0,\infty)\to[0,\infty)$ be a function with a compact support that is non-increasing in each variable. 
Then, the method we propose takes the form
\begin{align}
\VKL{} f:=
\sum_{k,\ell=0}^{\infty}
	\sum_{m=-\ell}^\ell
				h\left(\frac{k}{K},\frac{\ell}{N}\right)
				c_{k\ell m}J^{(\alpha,\beta)}_k Y_{\ell m}.
\end{align}
Note that this is actually a finite sum.
 Here, the coefficients $\{c_{k\ell m}\}\subset\mathbb{R}$ are quadrature approximations of Fourier coefficients,
\[
c_{k\ell m}
\approx
\frac1{\gamma_k^2}
\int_{\shell} f J^{(\alpha,\beta)}_k Y_{\ell m}.
\]
The quadrature approximation, the measure used in the integral, and the normalising constant $\gamma_k$ are defined later.

 Following \red{\cite{Sloan.I_Womersley_2012_filtered_hyperinterp,Wang.H_Sloan_2017_filtered}}, we shall call $\VKL{} f$ filtered hyperinterpolation of $f$ on $\shell$, if the quadrature is of suitably high polynomial precision.

Our main result Corollary \ref{cor:sup estim} gives error convergence orders of the method we propose in terms of the supremum norm. \red{Our strategy for the proof is similar to \cite{Sloan.I_Womersley_2012_filtered_hyperinterp,Wang.H_Sloan_2017_filtered} in that we reduce the error estimate to suitable best polynomial approximations. What differs in our setting is that we have the radial direction as well. We treat the error in the radial and angular directions separately by introducing} the filtered hyperinterpolation operator $\VK $ in the radial direction and $\VL$ in the angular direction. The error $f-\VKL{}f$ in terms of the supremum norm turns out to be bounded by the sum of the error bounds for each direction.

The outline of this paper is as follows. Section 2 introduces notations we need. In Section 3 and 4, we introduce the filtered hyperinterpolation approximations in the radial direction and the angular direction. Section 5 develops the filtered hyperinterpolation on spherical shells and analyses the error. We give numerical results in Section 6, and Section 7 concludes the paper.
\section{Preliminaries}
We set up some notations and introduce the problem we consider. 

With $\rin\in(0,1]$ and $\rout\in[1,\infty)$ ($\rin\neq\rout$), we consider a spherical shell $\shell=\{\bs{x}\in\mathbb{R}^3\ |\ \rin\le \norm{\bs{x}}_2\le\rout\}$. We assume $\rout-\rin=\ep>0$. We use the spherical coordinate system 
\[
\bs{x}\!=\!r\bs{\sigma}\!=(r\sin\theta\cos\varphi,r\sin\theta\sin\varphi,r\cos\theta)\quad(r\in[0,\infty),\ \! \theta\in[0,\pi],\ \varphi\in[0,2\pi)),\]
where $r=\norm{\bs{x}}_2$, $\bs{\sigma}=\frac{\bs{x}}{r}$, and for $\theta\in \{0,\pi\}$ we let $\varphi =0$.

For $\bs{\sigma}\in S^2$, we often write a function $f(\theta,\varphi)$ on the unit sphere $S^2$ as $f(\bs{\sigma})$.

In the following, we introduce orthogonal polynomials on the interval and the sphere. Further, we introduce the approximation method we consider.
\subsection{Orthogonal polynomials}
Let $J^*_k=J^{*(\alpha,\beta)}_k$ be the Jacobi polynomial of degree $k$ with the parameters $\alpha,\beta>-1$ on $[-1,1]$. Define $J_k$ $(k=0,\dotsc,{K})$ by
\[J_k(r)=J^*_k\left(\frac{2r-(r_{\rm in}+r_{\rm out})}{r_{\rm out}-r_{\rm in}}\right),\qquad r\in\rinout\]

 Let $w^*(x):=(1-x)^\alpha (1+x)^\beta$ ($x\in(-1,1)$) be the weight function associated with $\{J^*_k\}=\{J^{*(\alpha,\beta)}_k\}$, that is, with $\gamma_k=\left(\int_{-1}^1(J_k^*(x))^2 w^*(x)\mathrm{d}x\right)^{\frac12}$ we have
\begin{align}
	\int_{-1}^1 J_j^*(x)J_k^*(x)w^*(x){\mathrm d}x=\delta_{j k}\gamma_k^2,
\end{align}
where $\delta_{j k}=1$ if $j=k$ and $\delta_{j k}=0$ otherwise. For example, the weight associated with Legendre polynomials ($\alpha=\beta=0$) is $w^*_{\mathrm{Legendre}}(x)=1$ ($x\in(-1,1)$), and for Chebyshev polynomials ($\alpha=\beta=-\frac12$) we have  $w^*_{\mathrm{Chebyshev}}(x)=\frac{1}{\sqrt{1-x^2}}$ ($x\in(-1,1)$). We always consider a fixed pair of parameters $(\alpha,\beta)$, and thus we omit them in the notation $J^*_k$, and $J_k$.

%Define the weight function $w$ on $[\rin,\rout]$ as 
Define the measure $\mu^{\mathrm{rad}}$ on $\rinout$ by 
%
%\[w(r)=w^*\left(\frac{2r-(\rout+\rin)}{\rout-\rin}\right)\frac2{\rout-\rin}\frac1{r^2}.\]
\[\mu^{\mathrm{rad}}(A)
=
\int_A w^*\left(\frac{2r-(\rout+\rin)}{\rout-\rin}\right)\frac2{\rout-\rin}\dr,\]
for any Lebesgue measurable set $A$ in $\rinout$.
Then, we have 
\begin{align}
	\int_{\rin}^{\rout} J_j(r)J_k(r)\dmur^{\mathrm{rad}}(r)
%	&=
%	\int_{\rin}^{\rout} J_j(r)J_k(r)
%		w^*\left(\frac{2r-(\rout+\rin)}{\rout-\rin}\right)\frac2{\rout-\rin}\dr\\
	&=
	\int_{-1}^{1}J_j^*(x)J_k^*(x)w^*(x)\mathrm{d}x=\delta_{j k}\gamma_k^2.
\end{align}
Let $Y_{\ell,m}(\theta,\varphi)$ be the real spherical harmonics on the unit sphere $S^2$ defined by
\begin{align*}
\left\{ \begin{array}{ll}
Y_{\ell ,0}(\theta,\varphi)  &= \displaystyle\frac{1}{\sqrt{2\pi}}\mathcal{P}_{\ell , 0}(\cos\theta) \vspace{2pt}\\
Y_{\ell ,m}(\theta,\varphi)   &= \displaystyle\frac{1}{\sqrt{\pi}}\mathcal{P}_{\ell , m}(\cos\theta)\cos m \varphi \quad  (m=1,...,\ell) \vspace{2pt}\\
Y_{\ell ,-m}(\theta,\varphi) &=\displaystyle \frac{1}{\sqrt{\pi}}\mathcal{P}_{\ell , m}(\cos\theta)\sin m \varphi \quad  (m=1,...,\ell),
\end{array}\right.
\end{align*}
where $\mathcal{P}_{\ell , m}$ are defined as follows. Consider the associated Legendre polynomials
\begin{align*}
\left\{ \begin{array}{ll}
\displaystyle P_\ell^{0}(x) = P_\ell(x)  										   & (m=0) \vspace{3pt}\\
\displaystyle P_\ell^{m}(x) = (1-x^2)^{m/2}\ \frac{\mathrm{d}^m P_\ell}{\mathrm{d}x^m}(x)  & (m=1,...,\ell) \vspace{3.9pt}\\
\displaystyle P^{-m}_\ell(x) = (-1)^m \frac{(\ell-m)!}{(\ell+m)!} P^{m}_\ell(x)  & (m=1,...,\ell),
\end{array}\right.
\end{align*}
where $P_\ell(x)$ $(x\in[-1,1])$ is the Legendre polynomial of degree $\ell\in\{0\}\cup\mathbb{N}$. Then, $\{\mathcal{P}_{\ell,m}(x)\}_{\ell=0,1,2,\dotsc\atop m=-\ell,\dotsc,\ell}$ are defined by
\begin{align*}
\left\{\mathcal{P}_{\ell,m}(x)\equiv \sqrt{\frac{2\ell +1}{2}\frac{(\ell -m)!}{(\ell +m)!}}P^{m}_\ell (x) \right\} _{|m|\leq \ell,\ell =0,1,2,\cdots}.
\end{align*}

We often write the integral $\int_0^{2\pi}\!\!\!\int_0^\pi
f(\theta,\varphi)\sin\theta{\rm d}\theta{\rm d}\varphi$ as $\int_{S^2} f(\theta,\varphi)\dS$ or $\int_{S^2} f(\bs{\sigma})\dS$. The above normalisation gives us 
\begin{align}
\int_{S^2}
Y_{\ell,m}(\theta,\varphi)
Y_{\mu,\nu}(\theta,\varphi)\dS
=\delta_{l\mu}\delta_{m\nu}.
\end{align}
Finally, we let $\mathbb{P}_{k}([\rin,\rout])$, and $\mathbb{P}_{\ell}(S^2)$ be the space of polynomials of degree $\le k$ on $[\rin,\rout]$, and respectively the space of spherical polynomial of degree $\le\ell$ on $S^2$.
For details of orthogonal polynomials, see, for example, \cite{Szego.G_book_1975_fourth,Vilenkin.N_1978}.
\begin{comment}
Finally, we define the function space $\Ltmushell$ as the set of real, measurable functions satisfying
\begin{align}
\Ltmushell:
=
\left\{
f\from\mathbb{\shell}\to\mathbb{R}
\,\bigg|\,
\norm{f}_{\Ltmushell}^2:=
\intr\!\!\!\int_{S^2}
|f(r,\theta,\varphi)|^2
r)\dS<\infty
\right\}.
\end{align}
We consider the normed space $\Ltmushell/\{\norm{f}_{\Ltmushell}=0\}$. We use the same notation $\big(\Ltmushell,\norm{\cdot}_{\Ltmushell}\big)$  for this space of equivalence classes.
%
%
\end{comment}
%
Consider functions $h^{\mathrm{ang}}, h^{\mathrm{rad}}\from[0,\infty)\to[0,\infty)$ with $\supp{(h^{\mathrm{ang}})}\subset[0,a]$, and $\supp{(h^{\mathrm{rad}})}\subset[0,b]$ ($a,b\in(1,2]$). 
Further, we assume $h^{\mathrm{ang}}(s)=h^{\mathrm{rad}}(t)=1$ for $s,t\in[0,1]$.
Let us define the {\it filter} function $h\from[0,\infty)\times[0,\infty)\to[0,\infty)$ by
\begin{align}
h\from (s,t)\mapsto h^{\mathrm{rad}}(s)h^{\mathrm{ang}}(t).
\end{align}
We consider an approximation of a real-valued function $f$ on the shell $\rinout\times S^2$ of the form
\begin{align}\label{eq:LMNf}
\VKL{} f:=
\sum_{k,\ell=0}^\infty
	\sum_{m = -\ell }^\ell
	 h\left(\frac{k}{K},\frac{\ell}{L}\right)  c_{k\ell m}\,J_k(r)Y_{\ell,m}(\theta,\varphi),
\end{align}
where coefficients $\{c_{k\ell m}\}\subset\mathbb{R}$ shall be defined in Section \ref{sec:shell}, \eqref{eq:full scheme}. They are approximations of Fourier coefficients, that is,
\begin{align}\label{eq:c_klm}
	c_{k\ell m}
	&\approx
	\frac{1}{\gamma_k^2}\int_0^{2\pi}\!\!\!\!\int_0^\pi\!\!\! \intr\!\!
		f(r,\theta,\varphi)J_k(r)Y_{\ell,m}(\theta,\varphi)\dmur^{\mathrm{rad}}(r)\sin\theta\mathrm{d}\theta\mathrm{d}\varphi.
\end{align}
%
\begin{comment}
Consider a space of functions that satisfies $\int_0^{2\pi}\!\!\int_0^\pi\!\! \intr |f|^2 w(r)r^2\dr\sin\theta\mathrm{d}\theta\mathrm{d}\varphi<\infty$. In this space, the following bilinear form defines an inner product. 
%
\begin{align}
	\innprod{f}{g}_{L^2(\shell)}
	=\int_{S^2} \intr f(r,\theta,\varphi)\,g(r,\theta,\varphi)\dmur^{\mathrm{rad}}(r)\dS,
\end{align}
where the notation $\intr d\mu^r(r)=\intr w(r)r^2\dr$ and $\int_{S^2}\dS=\int_0^{2\pi}\!\!\int_0^\pi\sin\theta{\mathrm d}\theta{\mathrm d}\varphi$ are used. 

Our aim is to obtain a bound of the form
%
\begin{align}
\norm{f-\mathscr{P}_{ML}f}_{\Ltmu(\shell)}\le \mathscr{E}_1(M)+\mathscr{E}_2(L),
\end{align}
where $\mathscr{E}_1$ and $\mathscr{E}_2$ are decreasing functions with respect to $K$ and $L$, respectively. Our strategy is to introduce operators $\mathscr{P}_M$ and $\mathscr{P}_L$ that satisfy $\mathscr{P}_{ML}=\mathscr{P}_M\mathscr{P}_L$ as intermediate steps, and evaluate the right hand side of
%
\begin{align}\label{eq:decompose}
	\norm{f-\mathscr{P}_{ML}f}_{\Ltmu(\shell)}\le
	\norm{f-\mathscr{P}_{Nr}f}_{\Ltmu(\shell)}+
	\sup_{0\neq f\in\Ltmu(\shell)}\frac{\norm{\mathscr{P}_{Nr}f}_{\Ltmu(\shell)}}{\norm{f}_{\Ltmu(\shell)}}
	\norm{f-\mathscr{P}_{L}f}_{\Ltmu(\shell)}.
\end{align}
\end{comment}

\section{Filtered hyperinterpolation on the radial interval}
In this section we define filtered hyperinterpolation in the radial direction, and we will see that it is bounded as an operator from $C(\shell)$ to $C(\shell)$.

In order to develop properties of the filtered hyperinterpolation, as an intermediate step we define the {\it continuous} filtered approximation in the radial direction.

Let $\innprod{f}{g}_{\Ltmu(\rinout)}:=\intr fg\dmur^{\mathrm{rad}}$, and $\Ka=\max\{\lceil a K\rceil-1,K\}$ with $a\in(1,2]$.  For $f\in C([\rin,\rout])$, we define the filtered approximation $\VKconti f$ ($K\ge1$) by 
\begin{align}
\VKconti f(r)
&=
\sum_{k=0}^\infty \hrad\innprod{f}{\frac{J_k}{\gamma_k}}_{\Ltmu(\rinout)} \frac{J_k(r)}{\gamma_k}\label{eq:conti op rad}\\
&=
\sum_{k=0}^\Ka \hrad\innprod{f}{\frac{J_k}{\gamma_k}}_{\Ltmu(\rinout)} \frac{J_k(r)}{\gamma_k}\\
&=
\innprod{f}{ \mathcal{G}_K(\cdot,r) }_{\Ltmu(\rinout)},
\end{align}
where 
\begin{align}
\mathcal{G}_K(s,r):=
\sum\limits_{k=0}^{\Ka}\hrad
	\frac{J_k(s)}{\gamma_k} \frac{J_k(r)}{\gamma_k}
	\quad\text{ for }(s,r)\in\rinout\times\rinout.
\end{align}
Note that this is a finite sum.

The following proposition is standard.
\begin{prop}\label{prop:RK op norm}
For $\VKconti$ defined by \eqref{eq:conti op rad} we have 
\begin{align}
\|\VKconti\|_{C({\rinout})\to C({\rinout})}
=
\sup_{r\in\rinout}\intr |\mathcal{G}_K(s,r)|\dmur^{\mathrm{rad}}(s).\label{eq:conti op norm}
\end{align}
\end{prop}
\begin{proof}
%--------------------------------------------------------
%----------------- proof begins -------------------------
%--------------------------------------------------------
Clearly, we have $$\|\VKconti\|_{C({\rinout})\to C({\rinout})}
\le
\sup_{r\in\rinout}\intr |\mathcal{G}_K(s,r)|\dmur^{\mathrm{rad}}(s),$$ 
since
$
|\VKconti f(r)|
\le
\Big(\sup_{t\in\rinout} |f(t)| \Big)
\intr |\mathcal{G}_K(s,r)|\dmur^{\mathrm{rad}}(s)
$
for each $r\in\rinout$. Conversely, let
$$
\mathrm{sgn}_r(s)
:=
\begin{cases}
\frac{ \mathcal{G}_K(s,r) }{ |\mathcal{G}_K(s,r)| }&\mathcal{G}_K(s,r)\neq 0,\\
0																									 &\text{otherwise}.
\end{cases}
$$
\begin{sloppypar}
Since $\mathcal{G}_K(\cdot,r)$ is measurable, so is $\mathrm{sgn}_r(\cdot)$. 
Clearly, $\mathrm{sgn}_r(\cdot)\in L^1_\mu:=L^1_{\mu^{\mathrm{rad}}}(\rinout)$. 
Thus, there exists a sequence $\{g_n\}_n\subset C(\rinout)$ such that
$g_n\to \mathrm{sgn}_r$ $(n\to\infty)$ in $L^1_{\mu}$.
Therefore, there exists a subsequence $\{g_{n,k}\}_k\subset\{g_n\}$ such that
\end{sloppypar}
$$g_{n,k}(s)\to \mathrm{sgn}_r(s)\qquad (k\to\infty)\quad\text{ for $\mu^{\mathrm{rad}}$-almost every }s.$$
Replacing, if necessary, each function $g_{n,k}$ with 
$\min\{\, \max\{-1,g_{n,k}\},1\}$, we can assume $|g_{n,k}|\le 1$ independently of $k$, so that 
$$|g_{n,k}(\cdot)\mathcal{G}_K(\cdot,r)|\le \max\limits_{s,t\in\rinout} |\mathcal{G}_K(s,t)|.$$ 
Thus, from the dominated convergence theorem we have
\begin{align*}
\intr |\mathcal{G}_K(s,r)|\dmur^{\mathrm{rad}}(s)
&=
\intr \mathrm{sgn}_r(s)\mathcal{G}_K(s,r)\dmur^{\mathrm{rad}}(s) \\
&=
\lim_{k\to\infty}\intr g_{n,k}(s)\mathcal{G}_K(s,r)\dmur^{\mathrm{rad}}(s)\\
&=
\lim_{k\to\infty} \VKconti g_{n,k}(r)
\le
\|\VKconti\|_{C({\rinout})\to C({\rinout})}.
\end{align*}
%--------------------------------------------------------
%----------------- proof ends ---------------------------
%--------------------------------------------------------
\end{proof}
The following result essentially due to Mhaskar \cite{Mhaskar.H.N_2004_polynomial_op_finite_interval} establishes the uniform boundedness of $\{\VKconti\}_K$.
\begin{theorem}[\cite{Mhaskar.H.N_2004_polynomial_op_finite_interval}]\label{thm:Mhaskar}
Let the parameters $(\alpha,\beta)$ for the Jacobi polynomial satisfy $\alpha,\beta\ge-\frac12$, and let $\iota\ge \alpha+\beta+1$ be an integer. 
 Suppose that the filter function $h^{\mathrm{rad}}\from[0,\infty)\to[0,\infty)$ satisfies $h^{\mathrm{rad}}(s)=1$ for $s\in[0,1]$, $\mathrm{supp}(h^{\mathrm{rad}})\subset[0,a]$ for some $a\in(1,2]$.
Further, suppose that $h^{\mathrm{rad}}$ and its derivatives of all orders up to $\iota-1$ are absolutely continuous, and the $\iota$-th derivative is of bounded variation. Then, $\VKconti\from C({\rinout})\to C({\rinout})$ defined by \eqref{eq:conti op rad} satisfies
\begin{align}
\sup_{K\ge1}\|\VKconti\|_{C({\rinout})\to C({\rinout})}<\infty.
\end{align}
\end{theorem}
\begin{proof}
Following \cite[Proof of Theorem 3.1]{Mhaskar.H.N_2004_polynomial_op_finite_interval}, we see that the current assumption on the filter $h^{\mathrm{rad}}$ implies \cite[Theorem 3.1]{Mhaskar.H.N_2004_polynomial_op_finite_interval}, and thus from \cite[(2.10)]{Mhaskar.H.N_2004_polynomial_op_finite_interval} we have
$$
\sup_{r\in\rinout}\intr |\mathcal{G}_K(s,r)|\dmur^{\mathrm{rad}}(s)<\infty.
$$
Then, the statement follows from Proposition \ref{prop:RK op norm}.

\end{proof}
The filtered hyperinterpolation defined in the next section is obtained by approximating $\innprod{\cdot}{\cdot}_{\Ltmu(\rinout)}$ by the Gauss-type quadrature.
\subsection{Filtered hyperinterpolation}
We first recall the Gauss-Jacobi quadrature rule. Let $\kappa_Q:=\left\lceil\frac{ \Ka + K -1 }{2}\right\rceil$, and $\{r_j\}_{j=0,\dotsc,\kappa_Q}\subset [\rin,\rout]$ be the zeros of a Jacobi polynomial $J_{\kappa_Q+1}=J^{(\alpha,\beta)}_{\kappa_Q+1}$. Then, there exists the corresponding positive weight $\{w^{\mathrm{rad}}_j\}_{j=0,\dotsc,\kappa_Q}$ that defines a quadrature rule $\Qrad(2\kappa_Q+1)$ with the precision $2\kappa_Q+1(\ge \Ka + K)$. 
That is, for $f\in C(\rinout)$ let 
$\Qrad(2\kappa_Q+1)f:=
\sum_{j=0}^{\kappa_Q} f(r_j)\weightrad_j$. Then, we have
\begin{align}
\Qrad(2\kappa_Q+1)q&=
%\sum_{j=0}^{\kappa_Q} q(r_j)\weightrad_j = 
\intr q(r)w^*\left(\frac{2r-(\rout+\rin)}{\rout-\rin}\right)\frac2{\rout-\rin}\dr \notag\\
&=
\intr q(r)\mathrm{d}\mu^{\mathrm{rad}}(r)
 ,
\end{align}
{for any $q\in \mathbb{P}_{2\kappa_Q+1}([\rin,\rout])$, in particular, for any $q\in \mathbb{P}_{\Ka + K}([\rin,\rout])$.}

In the following, for a function $f\from\shell\to\bbR$ on $\shell$ we often write  $f(r\bs{\sigma})$ $(r\bs{\sigma}\in\shell)$ as $f(r,\bs{\sigma})$. For the sake of simplicity, we introduce the notation
\begin{align}
\innprod{f(\cdot,\bs{\sigma})}{g(\cdot,\bs{\sigma})}^{\mathrm{rad}}_{\Qrad(2\kappa_Q+1)}:=
%\innprod{f(\cdot,\bs{\sigma})}{g(\cdot,\bs{\sigma})}^{\mathrm{rad}}_{\Qrad(\kappa_Q+1;\Ka + K)}:=
\sum_{j=0}^{\kappa_Q} \weightrad_jf(r_j,\bs{\sigma})g(r_j,\bs{\sigma}),\ \ 
 f,g\in C(\shell),
\end{align}
where a quadrature rule $\Qrad(2\kappa_Q+1)$ is used on the right hand side.
\begin{sloppypar}
We defined the filtered approximation $\VKconti$ as \eqref{eq:conti op norm} with the calligraphic character. Here, we define the discretised filtered approximation operator $\VK \from C(\shell)\to C(\shell)$ for the radial direction, a quadrature-based approximation of $\VKconti$.
 For a function $f$ in $C(\shell)
,$
we define the operator $\VK $ ($K\ge1$) as
\end{sloppypar}
\begin{comment}
\begin{align}
\VK f(r,\bs{\sigma})
=
\sum_{k=0}^\infty \hrad\innprod{f(\cdot,\bs{\sigma})}{\frac{J_k}{\gamma_k}}^{\mathrm{rad}}_{\Qrad} \frac{J_k(r)}{\gamma_k}
\end{align}
%
Note that we have $\sum_{k=0}^\infty \hrad=\sum_{k=0}^{\Ka}\hrad+0$. Thus, we may define $\VK $ as
%
\end{comment}
\begin{align}
\VK f(r,\bs{\sigma})
=
\sum_{k=0}^\Ka \hrad\innprod{f(\cdot,\bs{\sigma})}{\frac{J_k}{\gamma_k}}^{\mathrm{rad}}_{\Qrad(2\kappa_Q+1)} \frac{J_k(r)}{\gamma_k}.\label{eq:def ang filtered}
\end{align}
Following \cite{Sloan.I_Womersley_2012_filtered_hyperinterp}, we call $\VK f$ {\it filtered hyperinterpolation} of $f$ in the radial direction. We now develop a bound of $\VK $ on the interval $\rinout$ in terms of the supremum norm over {$\rinout$}. Later, we use the result to analyse the error of filtered hyperinterpolation on spherical shells.
\subsection{Supremum norm bound of $\VK f$}
In this section, we obtain a bound of $\VK f$ in terms of the supremum norm. We have the following bound. 
\begin{prop}\label{prop:rad sup bd}
{Let the parameters $(\alpha,\beta)$ for the Jacobi polynomial satisfy $\alpha,\beta\ge-\frac12$, and let $\iota\ge \alpha+\beta+1$ be an integer. 
 Suppose that the filter function $h^{\mathrm{rad}}\from[0,\infty)\to[0,\infty)$ satisfies $h^{\mathrm{rad}}(s)=1$ for $s\in[0,1]$, $\mathrm{supp}(h^{\mathrm{rad}})\subset[0,a]$ for some $a\in(1,2]$.
Further, suppose that $h^{\mathrm{rad}}$ and its derivatives of all orders up to $\iota-1$ are absolutely continuous, and the $\iota$-th derivative is of bounded variation. }
For $f\in C(\shell)$, let $\VK f$ ($K\ge 1$) be defined by \eqref{eq:def ang filtered} with 
$\Ka=\max\{ \lceil aK\rceil -1,K \}$ ($a\in(1,2]$). 
Then, for each $\bs{\sigma}\in S^2$ we have
\begin{align}
\sup_{r\in\rinout}|\VK f(r,\bs{\sigma})|\le
C_1
\sup_{r\in\rinout}|f(r,\bs{\sigma})|,
\label{eq:point-wise bound VM f}
\end{align}
where the constant $C_1$ is independent of \red{$\bs{\sigma}$}, $K$ and $f$.
\end{prop}
\begin{proof}
Fix $\bs{\sigma}\in S^2$.
Clearly, we have
\begin{align}
|\VK f(r,\bs{\sigma})|
\le
(\sup_{r\in\rinout}|f(r,\bs{\sigma})|)
\left(
\sum_{j=0}^{\kappa_Q} \weightrad_j|\mathcal{G}_K(r_j,r)|
\right).\label{eq:VM f bound first}
\end{align}
Note that $\mathcal{G}_K(\cdot,r)$ is a polynomial of degree $\le \Ka$. From a non-trivial result by Nevai which gives a bound on {Gauss-Jacobi quadrature formulae in terms of the integral the quadrature approximates,} (see, for example \cite[p. 35, theorem 4.7.4]{Nevai.P_1986_case_study}) we have 
\begin{align}
\sum_{j=0}^{\kappa_Q} \weightrad_j|\mathcal{G}_K(r_j,r)|
\le
C\left(
\frac{\Ka}{\kappa_Q+1}+1
\right)
\intr |\mathcal{G}_K(s,r)|\dmur^{\mathrm{rad}}(s),
\end{align}
where $C$ depends only on the measure $\mu^{\mathrm{rad}}$.

Since 
$
\frac{\Ka +K -1}{2}\le
\kappa_Q
$, we have 
%$\frac{1}{(\kappa_Q+1)}\le \frac2{K+\Ka}$. Therefore, we have
\begin{align}
\frac{\Ka}{(\kappa_Q+1)}\le
\frac{2\Ka}{K+\Ka}\le 2,
\end{align}
and thus
\begin{align}
\sum_{j=0}^{\kappa_Q} \weightrad_j|\mathcal{G}_K(r_j,r)|
\le
3C
\intr |\mathcal{G}_K(s,r)|\dmur^{\mathrm{rad}}(s).
\end{align}
This, together with \eqref{eq:conti op norm} and \eqref{eq:VM f bound first}, we have
\begin{align}
\sup_{r\in\rinout}|&\VK f(r,\bs{\sigma})|\notag\\
&\le
3C\bigg(\sup_{r\in\rinout}|f(r,\bs{\sigma})|\bigg)
\sup_{r\in\rinout}
\intr |\mathcal{G}_K(s,r)|\dmur^{\mathrm{rad}}(s)\\
&=
3C\bigg(\sup_{r\in\rinout}|f(r,\bs{\sigma})|\bigg)
\|\VKconti\|_{C(\rinout)\to C(\rinout)}.\label{eq:VK bd by VKconti}
\end{align}
{In view of Theorem \ref{thm:Mhaskar},} we can conclude that 
\begin{align}
{\sup_{r\in\rinout}|\VK f(r,\bs{\sigma})|\le C_1\sup_{r\in\rinout}|f(r,\bs{\sigma})|,}
\end{align}
where the constant $C_1$ is independent of \red{$\bs{\sigma}$}, $K$ and $f$.
\end{proof}
%Note that, from the continuity of $f(\cdot,\bs{\sigma})$ on $\rinout$ we have $\sup_{r\in\rinout}|f(r,\bs{\sigma})|=\esssup_{r\in\rinout}|f(r,\bs{\sigma})|$, and thus 
%\begin{align}
%\sup_{r\in\rinout}|\VK f(r,\bs{\sigma})|\le c\esssup_{r\in\rinout}|f(r,\bs{\sigma})|\label{eq:esssup bound VM f}
%\end{align}
%holds.
\section{Filtered hyperinterpolation on sphere}
We now introduce filtered hyperinterpolation in the angular direction. Let $\Lb=\max\{\lceil bL\rceil-1,L\}$ with $b\in(1,2]$. Let $Q_{\mathrm{ang}}(\Lb + L)$ be a positive-weight $(\nu_Q+1)$-point spherical numerical integration rule
\begin{align}
	Q_{\mathrm{ang}}(\Lb + L) f(r):=
	\sumQN \weightang_j f(r,\bs{\sigma}_j),\quad f\in  C(\shell),
\end{align}
with points $\bs{\sigma}_0,\dotsc,\bs{\sigma}_{\nu_Q}\in S^2$ and corresponding weights $\weightang_0,\dotsc,\weightang_{\nu_Q}$, which integrates all spherical polynomials of degree $\le \Lb + L$ exactly. That is, we have
\begin{align}
 	\int_{S^2} p \dS=Q_{\mathrm{ang}}(\Lb + L) p,\quad p\in\mathbb{P}_{\Lb + L}(S^2).
\end{align}
Using this quadrature rule, let us define a bilinear map $\discipAng{\cdot}{\cdot}{\Lb + L}\from C(\shell) \times C(\shell) \to C(\shell)$ by
\begin{align}\label{eq:discinnprod}
\discipAng{f}{g}{\Lb + L}\!\!(r):=
%\left\langle f, g \right\rangle ^{\mathrm{ang}}_{Q_{\mathrm{ang}}(\nu_Q+1;L)}\!\!(r) :=
	\sumQN \weightang_j f(r,\bs{\sigma}_j)g(r,\bs{\sigma}_j),\quad f,g\in C(\shell).
\end{align}
Similarly to the radial direction, we consider a filter function $h^{\mathrm{ang}}\from[0,\infty)\to[0,\infty)$ with $h^{\mathrm{ang}}(s)=1$ for $s\in[0,1]$ and $\mathrm{supp}(h^{\mathrm{ang}})\subset[0,b]$ with a suitable smoothness. We define the filtered hyperinterpolation operator $\VL \from C(\shell)\to C(\shell)$ in the angular direction by 
\begin{align}\label{eq:defofLN}
	 (\VL f)(r,\bs{\sigma})
	 &:=
	 \sum_{\ell=0}^\infty\sum_{m=-\ell}^\ell
	 \hang
	 \discipAng{f(r,\cdot)}{Y_{\ell,m}(\cdot)}{\Lb+L}Y_{\ell,m}(\bs{\sigma})\\
	 &=
	 \sum_{\ell=0}^\Lb\sum_{m=-\ell}^\ell
	 \hang
	 \discipAng{f(r,\cdot)}{Y_{\ell,m}(\cdot)}{\Lb+L}Y_{\ell,m}(\bs{\sigma}).
\end{align}
%
\begin{comment}
Note that $\VN$ is idempotent. To see this, it suffice to check that for
\[g\in\mathscr{G}_N=\{g\in \Ltmushell:g(r,\theta,\varphi)=\sumLwithm  g^{\rm r}_{\ell,k}(r)Y_{\ell,m}(\bs{\sigma}),\,g^{\rm r}_{\ell,k}\in\Ltmu([\rin,\rout])\},\]
$\VN g=g$ holds. First, note that for each $\ell$ and $k$, we have
\begin{align}
\discipS{g(r,\cdot)}{Y_{\ell,m}}
=&
\sumQN\weightang_j
\Big(\sum_{\nu=0}^N\sum_{\mu=-\nu}^\nu
 g^{\rm r}_{\nu,\mu}(r)Y_{\nu,\mu}(\bs{\sigma}_j)
\Big)
Y_{\ell,m}(\bs{\sigma}_j)\\
=&
\sum_{\nu=0}^N\sum_{\mu=-\nu}^\nu
 g^{\rm r}_{\nu,\mu}(r)
\Big(
	\sumQN\weightang_j
	Y_{\nu,\mu}(\bs{\sigma}_j)
	Y_{\ell,m}(\bs{\sigma}_j)
\Big)\\
=&
\sum_{\nu=0}^N\sum_{\mu=-\nu}^\nu
 g^{\rm r}_{\nu,\mu}(r)
\Big(
	\int_{S^2}
	Y_{\nu,\mu}(\bs{\sigma})
	Y_{\ell,m}(\bs{\sigma})
	\dS
\Big)= g^{\rm r}_{\ell,k}(r).
\end{align}
%
Therefore,
%
\begin{align}
	\VN g(r,\bs{\sigma})
	=&\sumLwithm 
	\discipS{g(r,\cdot)}{Y_{\ell,m}}Y_{\ell,m}(\bs{\sigma})\\
	=&
	\sumLwithm 
	g^{\rm r}_{\ell,k}(r)Y_{\ell,m}(\bs{\sigma}).
\end{align}
Since $\VN f\in\mathscr{G}_N$ for $f\in D(L_N)$, we have $\VN(\VN f)=\VN f$.
\end{comment}
\subsection{Supremum norm estimate on the sphere}
To obtain our error estimate on spherical shells, we need the following supremum norm error estimate for the hyperinterpolation in the angular direction.
%------------------------------
%------------------------------
%----- POINT-WISE ESTIMATE ON SPHERE
%------------------------------
%------------------------------
\begin{theorem}[supremum norm estimate on $S^2$, %\cite{LeGia.T_Mhaskar_L_20082009_localized,Sloan.I_Womersley_2012_filtered_hyperinterp}]
%\cite{Mhaskar.H.N_2005_finitely_many,Sloan.I_Womersley_2012_filtered_hyperinterp}
\cite{Wang.H_Sloan_2017_filtered}]
\label{thm:sup estimate sphere} Let $r\in\rinout$. Suppose that the filter function $h^{\mathrm{ang}}\from[0,\infty)\to[0,\infty)$ satisfies $h^{\mathrm{ang}}(s)=1$ for $s\in[0,1]$, $\mathrm{supp}(h^{\mathrm{ang}})\subset[0,b]$ for some $b\in(1,2]$.
{Further, suppose that $h^{\mathrm{ang}}$ is absolutely continuous and its derivative is of bounded variation.} Then, for $f(r,\cdot)\in C(S^2)$ and $L\ge 1$, we have 
\begin{align}\label{eq:bd ang sup}
%\sup_{\bs{\sigma}\in S^2}|f(r,\bs{\sigma})-\VL f(r,\bs{\sigma})|\le c\|f(r,\cdot)\|_{W^s_\infty(S^2)} L^{-s},
\sup_{\bs{\sigma}\in S^2}|f(r,\bs{\sigma})-\VL f(r,\bs{\sigma})|\le
C_2
\inf_{p\in\mathbb{P}_L(S^2)}\sup_{\bs{\sigma}\in S^2}|f(r,\bs{\sigma})-p(\bs{\sigma})|.
\end{align}
where the constant $C_2$ {is independent of \red{$r$,} $f$ and $L$.}
%depends only on $h^{\mathrm{ang}}$, and $\iota$.
\end{theorem}
% See also \cite{LeGia.T_Mhaskar_L_20082009_localized}.
\begin{proof}
{From \cite[Proof of Theorem 1.1]{Wang.H_Sloan_2017_filtered},
for any $g\in C(S^2)$ we have 
$$
\sup_{\bs{\sigma}\in S^2}|\VL g(\bs{\sigma})|\le C \sup_{\bs{\sigma}\in S^2}| g(\bs{\sigma})|,$$
with a constant $C>0$ independent of $g$ and $L$. 
\red{In particular, for any $f(r,\cdot)\in C(S^2)$ we have 
$\sup_{\bs{\sigma}\in S^2}|\VL f(r,\bs{\sigma})|\le C \sup_{\bs{\sigma}\in S^2}| f(r,\bs{\sigma})|$ with the same constant $C$, which is independent of $r$, $f$, and $L$. 
}
Now, observe $\VL p(\bs{\sigma})=p(\bs{\sigma})$ for any $p\in\mathbb{P}_L(S^2)$. Then, the statement follows using the standard technique
 $$
 |f(r,\bs{\sigma})-\VL f(r,\bs{\sigma})|
 \le 
 |f(r,\bs{\sigma})-p(\bs{\sigma})|+|\VL(p(\bs{\sigma})-f(r,\bs{\sigma}))|,
 $$
 for arbitrary $p\in\mathbb{P}_L(S^2)$.}
\end{proof}
\section{Filtered hyperinterpolation on spherical shells}\label{sec:shell}
We finally define the filtered hyperinterpolation operator on $\shell$, and give an error estimate in terms of the supremum norm over $\shell$.

First, note that for $f\in C(\shell)$, we have $\VL\VK f=\VK \VL f$. We define the operator $\VKL{}$ on $C(\shell)$ as %D(\VK)\cap D(\VL)$ (i.e., separate continuous functions on $\shell$) as

\begin{align}
	(&\VKL{} f)(r,\bs{\sigma})\notag\\
	=&\!\sum_{k,\ell=0}^\infty
	\sum_{m=-\ell}^\ell
	h\bigg(\frac{k}{K},\frac{\ell}{L}\bigg)
	\Bigg(
	\sum_{j=0}^{\kappa_Q}
	\sum_{n=0}^{\nu_Q}
	W_{j n}f(r_j,\bs{\sigma}_n)
	\frac{J_k(r_j)}{\gamma_k^2}
	Y_{\ell,m}(\bs{\sigma}_n)
	\Bigg)
	{J_k(r)}Y_{\ell,m}(\bs{\sigma}),\label{eq:full scheme}
\end{align}
where $h\left(\frac{m}{M},\frac{\ell}{L}\right)=\hrad\hang$, $W_{j n}=\weightrad_j\weightang_n$, {$\kappa_Q=\Big\lceil\frac{\Ka+K-1}2\Big\rceil$} so that the radial quadrature has the precision $\Ka + K$, and $\nu_Q$ is taken so that the angular quadrature has the precision $L+\Lb$.

We estimate the error by following decomposition. For an arbitrary norm $\|\cdot\|$, from $\VKL{}=\VK \VL$ we have 
\begin{align}
 \norm{f-\VKL{} f}
 &=\norm{f-\VK \VL f+\VK f-\VK f}\nonumber\\
 &\le
 \norm{f-\VK f}
 +
 \norm{\VK (f-\VL f)}.
 \label{eq:decomperr}
\end{align}
%\blue{We do not want the other decomposition because we do not have the bound $\norm{\VL f}_{L^2\text{-sphere}}\le c\norm{f}_{C(\shell)}$.  We do not want to use $\norm{\VL f}_{L^2\text{-sphere}}\le c\norm{f}_{W^r(\shell)}$ as it makes the splitting complicated or it would lead us to use a $W^r(\shell)$ norm estimate}.

We derive estimates for both terms $\norm{f-\VK f}$ and $ \norm{\VK (f-\VL f)}$, with $\|\cdot\|$ being the supremum norm.
\subsection{Best approximation by polynomials}
We record classical results of estimates on best approximation by polynomials. Later, we reduce the error of the filtered hyperinterpolation approximation to the best approximation error.

On the interval, we have the following well-known results (see for example, \cite[pp. 196--197, p. 26]{Powell.J.D_book}).
\begin{theorem}\label{thm:best interval sup}
Let $$
\mathcal{E}^{\mathrm{rad}}_{K,\infty}(f)
=\inf_{p\in\mathbb{P}_K (\rinout)}\sup_{r\in\rinout}|f(r)-p(r)|.$$ Then, for $f\in C^\eta(\rinout)$ with $\eta\in\{1,\dotsc,K\}$, we have
\begin{align}
\mathcal{E}^{\mathrm{rad}}_{K,\infty}(f)\le C_3^\eta\frac{(K-\eta)!}{K!}
\bigg(\sup_{r\in\rinout}|f^{(\eta)}(r)|\bigg),
\end{align}
where the constant $C_3$ is independent of $f$ and $K$.
\end{theorem}
\begin{comment}
We record a version where the error is measured with respect to essential supremum. See, for example, \cite[(5.4.16)]{Canuto.C2006book} or \cite[Proof of Theorem 2]{Quarteroni.A_1984_Bernstein_Jackson_JJIAM}. 
\begin{theorem}\label{thm:best interval esssup}
Then, for $f\in W^\eta_{\infty}(\rinout)$ ({define somewhere}) with $1\le\eta\le M$, we have
\begin{align}
\inf_{p\in\mathbb{P}_M(\rinout)}\esssup_{r\in\rinout}|f(r)-p(r)|
\le CM^{-\eta}
\bigg(\max_{0\le j\le \eta}\esssup_{r\in\rinout}|D_r^jf(r)|\bigg),
\end{align}
where $D_r^jf(r)$ is the distributional derivative of order $j\in\{0,1\dotsc,\eta\}$ in $r$, and the constant $C$ is independent of $f$ and $K$.
\end{theorem}
(If $\esssup_{r\in\rinout}|D_r^\eta f(r)|<\infty$, then $|D_r^\eta f(r)|$ is square integrable on the set $(\rin,\rout)$ of finite measure. Thus $f\in W_2^\eta((\rin,\rout))$. Thus $f^{(\eta-1)}$ is absolutely continuous.)
\end{comment}

On the sphere, we have the following classical result by Pawelke \cite{Pawelke_1972}. 
\begin{theorem}[Pawelke \cite{Pawelke_1972}]\label{thm:best sphere}
Let $$
\mathcal{E}^{\mathrm{ang}}_{L,\infty}(f)
=
\inf_{p\in\mathbb{P}_L(S^2)}\sup_{\bs{\sigma}\in S^2}|f(\bs{\sigma})-p(\bs{\sigma})|
.
$$
Then, for each $f\in C^{2t}(S^2)$, ($t=1,2,\dotsc$), there exists a constant $C_4>0$ independent
of $f$ and $L$, such that
\begin{align}
\mathcal{E}^{\mathrm{ang}}_{L,\infty}(f)
\le C_4^tL^{-2t}\bigg(\sup_{\bs{\sigma}\in S^2}|\Delta_S^tf(\bs{\sigma})|
	\bigg)
\end{align}
holds, where $\Delta_S$ is the Laplace--Beltrami operator on $S^2$.
\end{theorem}
%{to be checked}Further, see Mhaskar, H. N., Narcowich, F. J., Prestin, J. and Ward, J. D. (2010). $L_p$ Bernstein estimates and approximation by spherical basis functions. Math. Comp. p. 1662, 

%Kamzolov, A. I. (1982). The best approximation of classes of functions $W^\alpha$ by polynomials in spherical harmonics. Mat. Zametki 32 285–293, 425.
%\begin{theorem}\label{thm:best sphere esssup}
%For each $f\in C^{2t}(S^2)$, ($t=1,2,\dotsc$), there exists a constant $C$ independent
%of $f$ and $L$, such that
%\begin{align}
%\inf_{p\in\mathbb{P}_L(S^2)}\sup_{\bs{\sigma}\in S^2}|f(r,\cdot)-p(\bs{\sigma})|\le C^tL^{-2t}\bigg(\sup_{\bs{\sigma}\in S^2}|\Delta_S^tf(\bs{\sigma})|
%	\bigg)
%\end{align}
%holds.
%\end{theorem}
\begin{rem}
\red{Note that we could also use the recent result on the best polynomial approximation on the sphere by Dai and Xu \cite[Corollary 3.7]{Dai.F_2011_Sobolev}, see also \cite[Corollary 4.5.6]{Dai.F2013book}. They considered a filtered approximation with a smooth filter, and reduced the error estimate to the best polynomial approximation. To analyse the convergence rate of the best polynomial approximation, they introduced a new class of Sobolev spaces on the sphere. In this paper, we adopt the classical result as above.}
\end{rem}

\subsection{Error estimate}
We have the following estimate.
\begin{theorem}
%Suppose that the quadrature rule used for $\VKL{} f$ is of degree of precision at least $\Ka+K$ in the radial direction, and in the angular direction $\Lb+L$. 
For $f\in C(\shell)$, let $\VKL{} f$ be defined by \eqref{eq:full scheme}. Then, under the same assumptions as Proposition \ref{prop:rad sup bd} and Theorem \ref{thm:sup estimate sphere}, we have
\begin{align}
\sup\limits_{(r,\bs{\sigma})\in\shell}&|f(r,\bs{\sigma})-\VKL{} f(r,\bs{\sigma})|\notag\\
&\le
(1+C_1)\sup_{\bs{\sigma}\in S^2}
\mathcal{E}^{\mathrm{rad}}_{K,\infty}(f(\cdot,\bs{\sigma}))
+
C_1C_{\red{2}}\sup_{r\in\rinout}\mathcal{E}^{\mathrm{ang}}_{L,\infty}(f(r,\cdot)),
\label{eq:bd best approx}
\end{align}
for $f\in C(\shell)$, where the constants $C_1$, $C_{\red{2}}>0$ are independent of \red{$r$, $\bs{\sigma}$}, $K$, $L$, and $f$.
\begin{comment}
Further, we have
\begin{align}
\sup\limits_{(r,\bs{\sigma})\in\shell}|f(r,\bs{\sigma})-\VKL{} f(r,\bs{\sigma})|
\le
(1+c)\esssup_{\bs{\sigma}\in S^2}
\mathcal{E}^{\mathrm{rad}}_{K,\infty}(f(\cdot,\bs{\sigma}))
+
c(1+c')\esssup_{r\in\rinout}\mathcal{E}^{\mathrm{ang}}_{L,\infty}(f(r,\cdot)).
\label{eq:bd best approx esssup}
\end{align}
\end{comment}
\end{theorem}
\begin{proof}
We first obtain a bound for the first term of \eqref{eq:decomperr} with the supremum norm over $\shell$. For an arbitrary $p\in\PolyRad{K}$, we have
\begin{align}
\sup_{(r,\bs{\sigma})\in\shell}|f(r,\bs{\sigma})-\VK f(r,\bs{\sigma})|
%&\blue{=\sup_{
%\substack{r\in\rinout\\\bs{\sigma}\in S^2}
%}|f(r,\bs{\sigma})-\VK f(r,\bs{\sigma})|}\\
%&\blue{=
%\sup_{\bs{\sigma}\in S^2}\sup_{r\in\rinout}
%	|f(r,\bs{\sigma})-\VK f(r,\bs{\sigma})|}\\
&\le
\sup_{\bs{\sigma}\in S^2}\sup_{r\in\rinout}
	(|(f-p)| + |\VK (f-p)|).
\end{align}
From \eqref{eq:point-wise bound VM f}, it follows that
\begin{align}
\sup_{\bs{\sigma}\in S^2}\sup_{r\in\rinout}
	(|(f-p)|& + |\VK (f-p)|)\notag\\
&\le
(1+C_1)\sup_{\bs{\sigma}\in S^2}
\sup_{r\in\rinout}|f(r,\bs{\sigma})-p(r)|.
\end{align}
Therefore, we have
\begin{align}
\sup_{(r,\bs{\sigma})\in\shell}|f(r,\bs{\sigma})-\VK f(r,\bs{\sigma})|
\le
(1+C_1)\sup_{\bs{\sigma}\in S^2}
\mathcal{E}^{\mathrm{rad}}_{K,\infty}(f(\cdot,\bs{\sigma})).\label{eq:bound 1st}
\end{align}
%
\begin{comment}
Note that we also have
\begin{align}
\sup_{(r,\bs{\sigma})\in\shell}|f(r,\bs{\sigma})-\VK f(r,\bs{\sigma})|
\le
(1+c)\esssup_{\bs{\sigma}\in S^2}
\mathcal{E}^{\mathrm{rad}}_{K,\infty}(f(\cdot,\bs{\sigma})).\label{eq:bound 1st esssup}
\end{align}
\blue{
Indeed, since $f(r,\bs{\sigma})-\VK f(r,\bs{\sigma})\in C(\shell)$ we have
\begin{align}
\sup_{(r,\bs{\sigma})\in\shell}|f(r,\bs{\sigma})-\VK f(r,\bs{\sigma})|
&=\esssup_{(r,\bs{\sigma})\in\shell}|f(r,\bs{\sigma})-\VK f(r,\bs{\sigma})|\\
&=\esssup_{
\substack{r\in\rinout\\\bs{\sigma}\in S^2}
}|f(r,\bs{\sigma})-\VK f(r,\bs{\sigma})|\\
&=
\esssup_{\bs{\sigma}\in S^2}\esssup_{r\in\rinout}
	|f(r,\bs{\sigma})-\VK f(r,\bs{\sigma})|\\
&\le
\esssup_{\bs{\sigma}\in S^2}\esssup_{r\in\rinout}
	(|\VK (f-p)|+|(f-p)|)\\
&\le
\esssup_{\bs{\sigma}\in S^2}\sup_{r\in\rinout}
	(|\VK (f-p)|+|(f-p)|).
\end{align}
Then, \eqref{eq:bound 1st esssup} follows from \eqref{eq:point-wise bound VM f}. %{maybe we do not even need \eqref{eq:esssup bound VM f} and make the inside esssup sup}
}
\end{comment}
On the other hand, again from \eqref{eq:point-wise bound VM f} we have

\begin{align}
\sup_{(r,\bs{\sigma})\in\shell}\big|\VK &\big(f(r,\bs{\sigma})-\VL f(r,\bs{\sigma})\big)\big| \notag\\
%
%&\blue{=
%\sup_{\bs{\sigma}\in S^2}\sup_{r\in\rinout}
%\big|\VK \big(f(r,\bs{\sigma})-\VL f(r,\bs{\sigma})\big)\big|\notag}\\
%
&\le
C_1 \sup_{\bs{\sigma}\in S^2}\sup_{r\in\rinout}
|f(r,\bs{\sigma})-\VL f(r,\bs{\sigma})|\label{eq:lhs 2nd term}\\
&=
C_1 \sup_{r\in\rinout}\sup_{\bs{\sigma}\in S^2}
|f(r,\bs{\sigma})-\VL f(r,\bs{\sigma})|.
\end{align}
\red{In view of Theorem 4.1, it follows that}
\begin{align}
C_1\sup_{r\in\rinout}\sup_{\bs{\sigma}\in S^2}
|f(r,\bs{\sigma})-\VL f(r,\bs{\sigma})|
\le
C_1\red{C_2}\sup_{r\in\rinout}\mathcal{E}^{\mathrm{ang}}_{L,\infty}(f(r,\cdot)).
\end{align}
Thus, the left hand side of \eqref{eq:lhs 2nd term} can be bounded as 
\begin{align}
\sup_{(r,\bs{\sigma})\in\shell}\big|\VK \big(f(r,\bs{\sigma})-\VL f(r,\bs{\sigma})\big)\big|
\le
C_1\red{C_2}\sup_{r\in\rinout}\mathcal{E}^{\mathrm{ang}}_{L,\infty}(f(r,\cdot)).\label{eq:bound 2nd}
\end{align}
%
\begin{comment}
We also have
\begin{align}
\sup_{(r,\bs{\sigma})\in\shell}\big|\VK \big(f(r,\bs{\sigma})-\VL f(r,\bs{\sigma})\big)\big|
\le
c(1+c')\esssup_{r\in\rinout}\mathcal{E}^{\mathrm{ang}}_{L,\infty}(f(r,\cdot)).\label{eq:bound 2nd esssup}
\end{align}
Indeed, we have
\begin{align}
\sup_{(r,\bs{\sigma})\in\shell}\big|\VK \big(f(r,\bs{\sigma})-\VL f(r,\bs{\sigma})\big)\big|
%
&=
\esssup_{(r,\bs{\sigma})\in\shell}\big|\VK \big(f(r,\bs{\sigma})-\VL f(r,\bs{\sigma})\big)\big|\\
&\le\esssup_{\bs{\sigma}\in S^2}\sup_{r\in\rinout}
\big|\VK \big(f(r,\bs{\sigma})-\VL f(r,\bs{\sigma})\big)\big|\\
%
&\le
c\esssup_{\bs{\sigma}\in S^2}\sup_{r\in\rinout}
|f(r,\bs{\sigma})-\VL f(r,\bs{\sigma})|.
\end{align}
From the continuity of $|f(r,\bs{\sigma})-\VL f(r,\bs{\sigma})|$, it follows that
\begin{align}
\esssup_{\bs{\sigma}\in S^2}\sup_{r\in\rinout}
|f(r,\bs{\sigma})-\VL f(r,\bs{\sigma})|
&=
\esssup_{\bs{\sigma}\in S^2}\esssup_{r\in\rinout}
|f(r,\bs{\sigma})-\VL f(r,\bs{\sigma})|\\
&=
\esssup_{r\in\rinout}\esssup_{\bs{\sigma}\in S^2}
|f(r,\bs{\sigma})-\VL f(r,\bs{\sigma})|\\
&\le
\esssup_{r\in\rinout}\sup_{\bs{\sigma}\in S^2}
|f(r,\bs{\sigma})-\VL f(r,\bs{\sigma})|,
\end{align}
and thus \eqref{eq:bound 2nd esssup} holds.
\end{comment}
%
From \eqref{eq:decomperr}, \eqref{eq:bound 1st}, and \eqref{eq:bound 2nd} the result \eqref{eq:bd best approx} follows. %The second claim \eqref{eq:bd best approx esssup} can be shown similarly.

\end{proof}
Together with Theorem \ref{thm:best interval sup} and \ref{thm:best sphere}, we have the following corollary of the previous theorem.
\begin{cor}\label{cor:sup estim}
For $f\in C(\shell)$, let $\VKL{} f$ be defined by \eqref{eq:full scheme}. Suppose the same assumptions as Proposition \ref{prop:rad sup bd} and Theorem \ref{thm:sup estimate sphere} hold. 
Further, suppose that $f$ is $\eta$-times continuously partially differentiable with respect to $r$ ($\eta\in\{1,\dotsc,K\}$) and satisfies $$\|f^{(\eta,0)}\|_{\infty,\infty}:=
\sup_{\bs{\sigma}\in S^2}\sup_{r\in\rinout}
\bigg|\frac{\partial^\eta}{\partial r^\eta}f(r,\bs{\sigma})\bigg|<\infty.$$
Suppose furthermore that $f(r,\cdot)\in C^{2t}(S^2)$ ($t\in\{1,2,3,\dots\}$) for each $r$ and satisfies
$$\|f^{(0,2t)}\|_{\infty,\infty}:=\sup_{r\in\rinout}\sup_{\bs{\sigma}\in S^2}|\Delta^t_Sf(r,\bs{\sigma})|<\infty.$$
Then, we have
\begin{align}
\sup\limits_{(r,\bs{\sigma})\in\shell}&|f(r,\bs{\sigma})-\VKL{} f(r,\bs{\sigma})|\notag\\
& \le
(1+C_1)C_3^\eta\frac{(K-\eta)!}{K!}\|f^{(\eta,0)}\|_{\infty,\infty}
+
C_1 \red{C_2}
\,C_4^tL^{-2t}\|f^{(0,2t)}\|_{\infty,\infty}
,
\end{align}
\red{where the constants $C_1,C_2,C_3,C_4>0$ are independent of $r$, $\bs{\sigma}$, $K$, $L$, and $f$.}
\begin{comment}
2) Suppose $f$ is $\eta$-times continuously partial differentiable with respect to $r$ ($\eta\in\{1,\dotsc,K\}$) and satisfies $$\|f^{(\eta,0)}\|_{\infty,\esssup}:=
\esssup_{\bs{\sigma}\in S^2}\sup_{r\in\rinout}
\bigg|\frac{\partial^\eta}{\partial r^\eta}f(r,\bs{\sigma})\bigg|<\infty.$$
Further, suppose $f(r,\cdot)\in C^{2t}(S^2)$ ($t\in\{1,2,3,\dots\}$) for each $r$ and satisfies
%
$$\|f^{(0,2t)}\|_{\esssup,\infty}:=\esssup_{r\in\rinout}\sup_{\bs{\sigma}\in S^2}|\Delta^t_Sf(r,\bs{\sigma})|<\infty.$$
Then, we have
\begin{align}
\sup\limits_{(r,\bs{\sigma})\in\shell}|f(r,\bs{\sigma})-\VKL{} f(r,\bs{\sigma})|
 \le
(1+c)c^\eta\frac{(K-\eta)!}{K!}\|f^{(\eta,0)}\|_{\infty,\esssup}
+
c(1+c')
\,C^tL^{-2t}\|f^{(0,2t)}\|_{\esssup,\infty}
.
\end{align}
\end{comment}
\end{cor}
\subsection{\red{Comparison with the non-filtered approximation in \cite{Kazashi.Y_2016_nonfiltered}}}
\label{sec:comparison}
\red{In \cite{Kazashi.Y_2016_nonfiltered} the author considered a fully discretised polynomial approximation on the shells; polynomial interpolation in the radial, and hyperinterpolation in the angular direction. 
We note that, when approximating smooth functions, the non-filtered approximation considered in \cite{Kazashi.Y_2016_nonfiltered} is not substantially worse than the filtered hyperinterpolation considered in this paper.}

\red{
We can derive an error estimate in terms of the supremum norm for the method considered in \cite{Kazashi.Y_2016_nonfiltered} following the argument in this paper. 
Then, in comparison with the method proposed in this paper, the convergence rate in terms the supremum norm will get worse only up to the factor of the product of the operator norms for the radial and angular approximations---interpolation and hyperinterpolation---as an operator from $C(\shell)$ to $C(\shell)$. 
The affinely mapped Chebyshev zeros as the interpolation points for the radial direction, for example, will yield the bound $O(\log K)$ for the operator norm in the radial direction, where $K$ is the highest degree of polynomial used for the radial direction. For the angular direction, under a mild condition on the spherical quadrature points, we have the bound $O(L^{\frac12})$ (see \cite[Theorem 5.5.4]{Sloan.I.H_Womersley_2000_constructive}), where $L$ is the highest degree of spherical harmonics used for the angular direction. For smooth functions, the effect of these factors will be insignificant relative to the convergence rate of the best polynomial approximations.}

\red{On the other hand, as the numerical results in the next section demonstrate, for non-smooth functions the filtered approximation works better.
 }
\section{Numerical results}
In this section, we provide numerical results. Let $[\rin,\rout]=[1,1.001]$ so that $\ep=0.001$.  
We let $a=b=2$ so that $\overline{K}(2)=2K-1$ ($K\ge1$), and $\overline{L}(2)=2L-1$ ($L\ge 1$). 
For the Jacobi polynomials we use Chebyshev polynomials of the first kind, that is, $(\alpha,\beta)=(-\frac12,-\frac12)$. For the radial direction, we use Gauss-Chebyshev quadrature with $\left\lceil \frac{3K-1}2 \right\rceil+1$ points $x_{0},\dotsc,x_{\left\lceil \frac{3K-1}2 \right\rceil}$ mapped to $\rinout$, and weights $w^x_{0},\dotsc,w^x_{\left\lceil \frac{3K-1}2 \right\rceil}$. For the angular direction, we use the {\it Spherical $t$-designs with $N$ = $t^2/2 + t + O(1)$ points} given by Womersley \cite{Womersley.R_2016_sph_des_web}. As for the filter, we use the following $C^\infty$ exponential filter proposed in \cite{Filbir.F_MP_2009_filter} for both directions:
\begin{align}
h^{\mathrm{rad}}(x),h^{\mathrm{ang}}(x)
:=
\begin{cases}
1 & \text{for }x\in [0,1]\\
\exp\Big(
	\frac
	{ -2\exp\big(\frac2{1-x}\big) }
	{ 2 - x }
	\Big) &\text{for }x\in (1,2)\\
0&\text{for }x\in [2,\infty).
\end{cases}
\end{align}
\red{For a comparison, we provide the error plot using the discretised polynomial approximation on spherical shells proposed in \cite{Kazashi.Y_2016_nonfiltered}. Specifically, for the comparison we employ interpolation at the Chebyshev zeros in the radial direction, and hyperinterpolation on the sphere in the angular direction. The error using this method will be simply referred to as (non-filtered) hyperinterpolation error.}

Let $r\in[1,1.001]$, $\theta\in[0,\pi]$, and $\varphi\in[0,2\pi)$. 
We first approximate $f_1(r,\theta,\varphi)=\frac{1}{0.0005^{\frac{17}2}}|r-1.0005|^{\frac{17}2}\cos^4\theta$, and 
$f_2(r,\theta,\varphi)=\frac{1}{0.001^2}(r-1)^{2}|\cos\theta|^{\frac{17}2}$. The function $f_1$ is smooth in the angular but less so in the radial direction; $f_2$ is smooth in the radial but less so in the angular direction.
\begin{figure}[H]%[htbp]%[H]
  \centering
  \includegraphics[width=\textwidth]{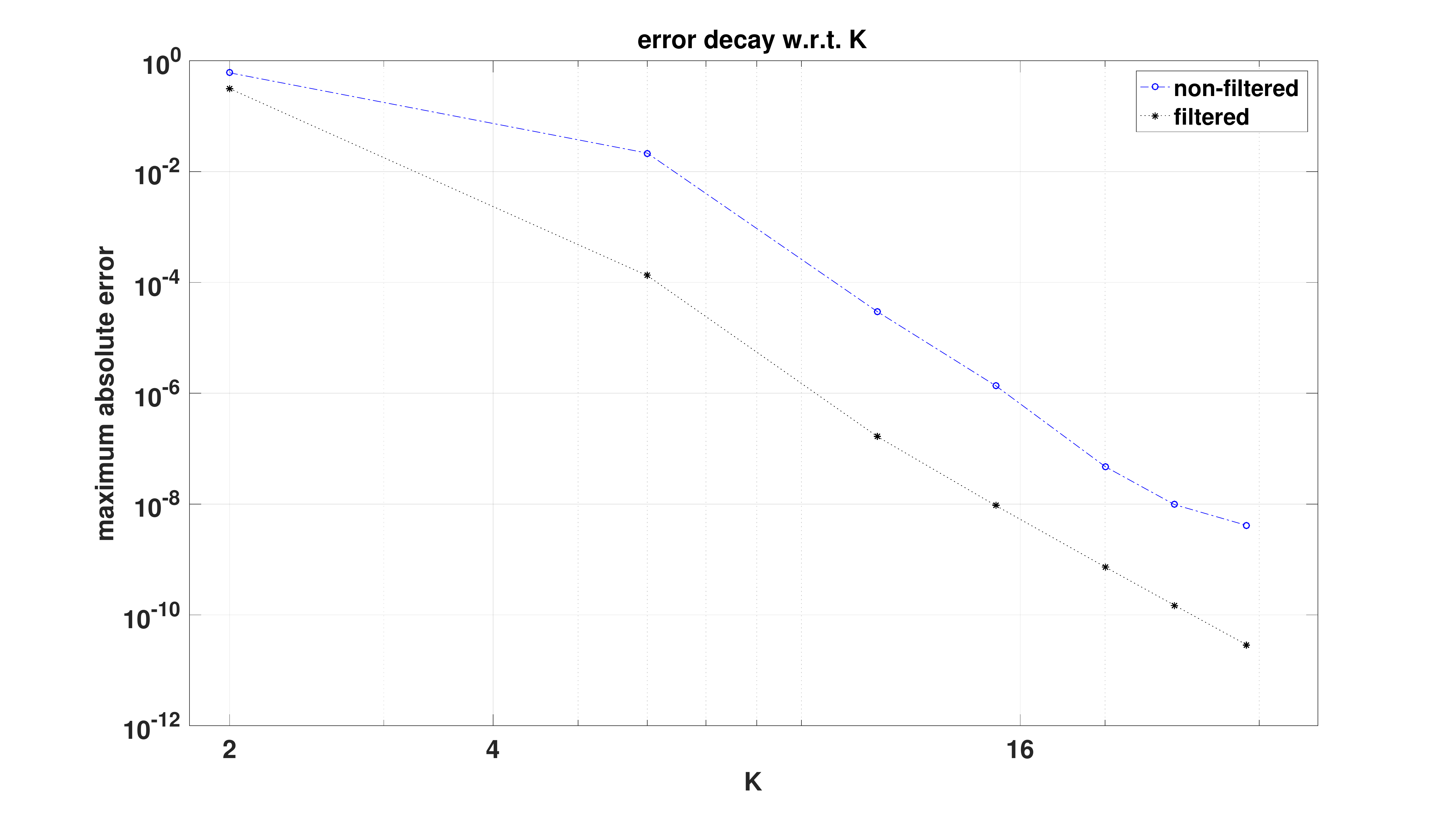}
  \caption{Approximand is $f_1(r,\theta,\varphi)=\frac{1}{0.0005^{\frac{17}2}}|r-1.0005|^{\frac{17}2}\cos^4\theta$. We fix $L=4$ and change $K$. We observe the algebraic decay of the error with respect to $K$ \red{for both filtered and non-filtered hyperinterpolations}.}
  \label{fig:change K}
\end{figure}
To approximate the function $f_1$, we fix $L=4$ and vary the value of $K$. We observe the algebraic decay of the maximum absolute error with respect to $K$ \red{for both filtered and non-filtered hyperinterpolations}. See Fig.~\ref{fig:change K}. The function $f_1$ is a $C^8(\rinout)$-function for a fixed angular variable. 

\red{As discussed in Section \ref{sec:comparison}, for this reasonably smooth function $f_1$ we observe that non-filtered hyperinterpolation converges at a rate not much worse than the filtered case. 
For the filtered approximation, we observe an almost the same convergence rate, but with a smaller error. 
Since $f_1$ is in $C^8(\rinout)$}, from Corollary \ref{cor:sup estim} we expect that the logarithm of the error decays no slower than $-8\log_{10}K$. The experiment shows roughly $\log_{10}(\mathrm{error})\approx-10\log_{10}K$.
 This seems to show our result might not be sharp, but at least it is consistent with the experiment.

\begin{figure}[H] %[htbp]%[H]
  \centering
  \includegraphics[width=\textwidth]{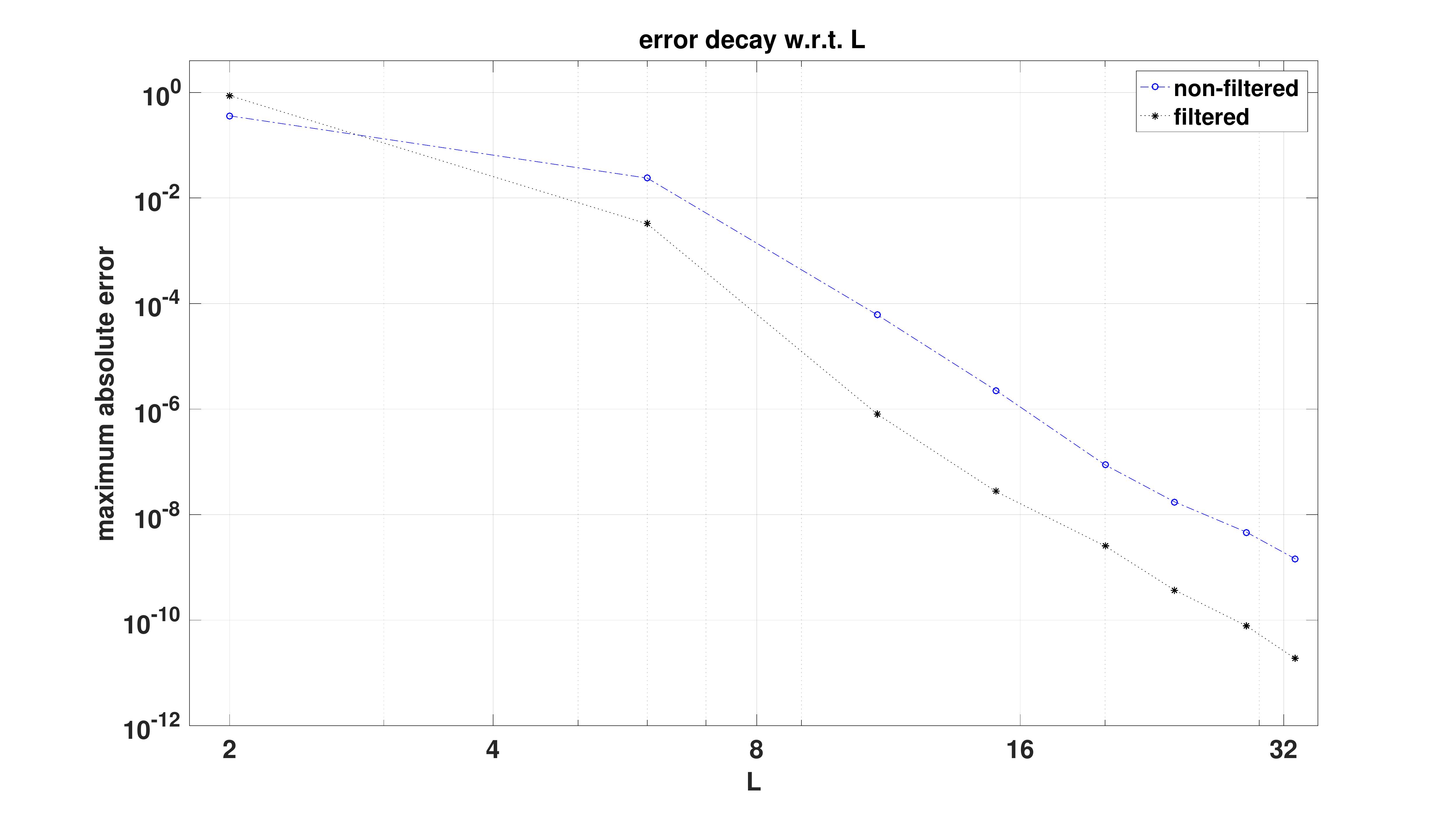}
  \caption{Approximand is $f_2(r,\theta,\varphi)=\frac{1}{0.001^2}(r-1)^{2}|\cos\theta|^{\frac{17}2}$. We fix $K=2$ and change $L$. We observe algebraic decay of the maximum absolute error with respect to $L$ \red{for both filtered and non-filtered hyperinterpolations}.}
  \label{fig:change L}
\end{figure}

Next, we approximate the function $f_2$. We fix $K=2$, so that the radial part can be approximated exactly, and vary the value of $L$. See Fig.~\ref{fig:change L}. We observe the algebraic decay of the error with respect to $K$. \red{The function $f_2$ is a $C^8(S^2)$-function for a fixed radial variable.}

\red{As discussed in Section \ref{sec:comparison}, for this reasonably smooth function $f_2$ we observe that non-filtered hyperinterpolation converges at a rate not much worse than the filtered case.
For the filtered approximation, we observe an almost the same convergence rate, but with a smaller error.} 
\red{Since $f_2$ is in $C^8(S^2)$}, from Corollary \ref{cor:sup estim} we expect the logarithm of the error decays no slower than $-8\log_{10}L$. The experiment shows that roughly $\log_{10}(\mathrm{error})\approx-11\log_{10}L$, which is consistent with Corollary \ref{cor:sup estim}.
Next, we consider a non-smooth function which goes beyond the theory. It is constructed as a sum of following two functions. First, the Franke function $f_{\mathrm{Franke}}$ is a $C^\infty$-function defined by
\begin{align}
f_{\mathrm{Franke}}(\bs{x})=&
0.75\exp\left(-\frac14((9x_1-2)^2 + (9x_2-2)^2 + (9x_3-2)^2)\right)\notag\\
& +        0.75\exp\left(-\frac1{49}((9x_1+1)^2) - \frac1{10}(9x_2+1+9x_3+1)\right)\notag\\
& + 0.5\exp\left(-\frac14((9x_1-7)^2 + (9x_2-3)^2 + (9x_3-5)^2)\right)\notag\\
& - 
         0.2\exp\left(-(9x_1-4)^2 - (9x_2-7)^2 - (9x_3-5)^2\right), 
\end{align}
for $\bs{x}=(x_1,x_2,x_3)\in\shell$.
Further, we define the cone function $f_{\mathrm{cone}}$ as follows
\begin{align}
f_{\mathrm{cone}}(\bs{x})
:=
\begin{cases}
1000\big|\|\bs{x}\|_2-1.0005\big|\Big(1-2\arccos\Big(\frac{\bs{x}}{\|\bs{x}\|_2}\cdot\bs{x}_{\mathrm{c}}\Big)\Big)
&\bs{x}\in \mathcal{C}(\bs{x}_{\mathrm{c}})\\
0 & \text{otherwise},
\end{cases}
\end{align}
\begin{sloppypar}
where $\mathcal{C}(\bs{x}_{\mathrm{c}}):=\{\bs{x}\in\shell\mid \arccos(\bs{x}/\|\bs{x}\|_2\cdot\bs{x}_{\mathrm{c}}))\le \frac12\}$, and $\bs{x}_{\mathrm{c}}:=(-1/2,-1/2,1/\sqrt{2})\in S^2$. We approximate 
\end{sloppypar}
\begin{align}
f_3(\bs{x}):=f_{\mathrm{Franke}}(\bs{x})+f_{\mathrm{cone}}(\bs{x}).
\end{align}
Note that $f_3$ is not differentiable in the angular direction along the boundary of the spherical cap $\mathcal{C}(\bs{x}_{\mathrm{c}})$ nor at the centre of the cap, and in the radial direction at the midpoint of the interval $[1,1.001]$.
\begin{figure}[H]%[htbp]%[H]
\centering
\subfloat[]{
   \includegraphics[width=0.50\textwidth]{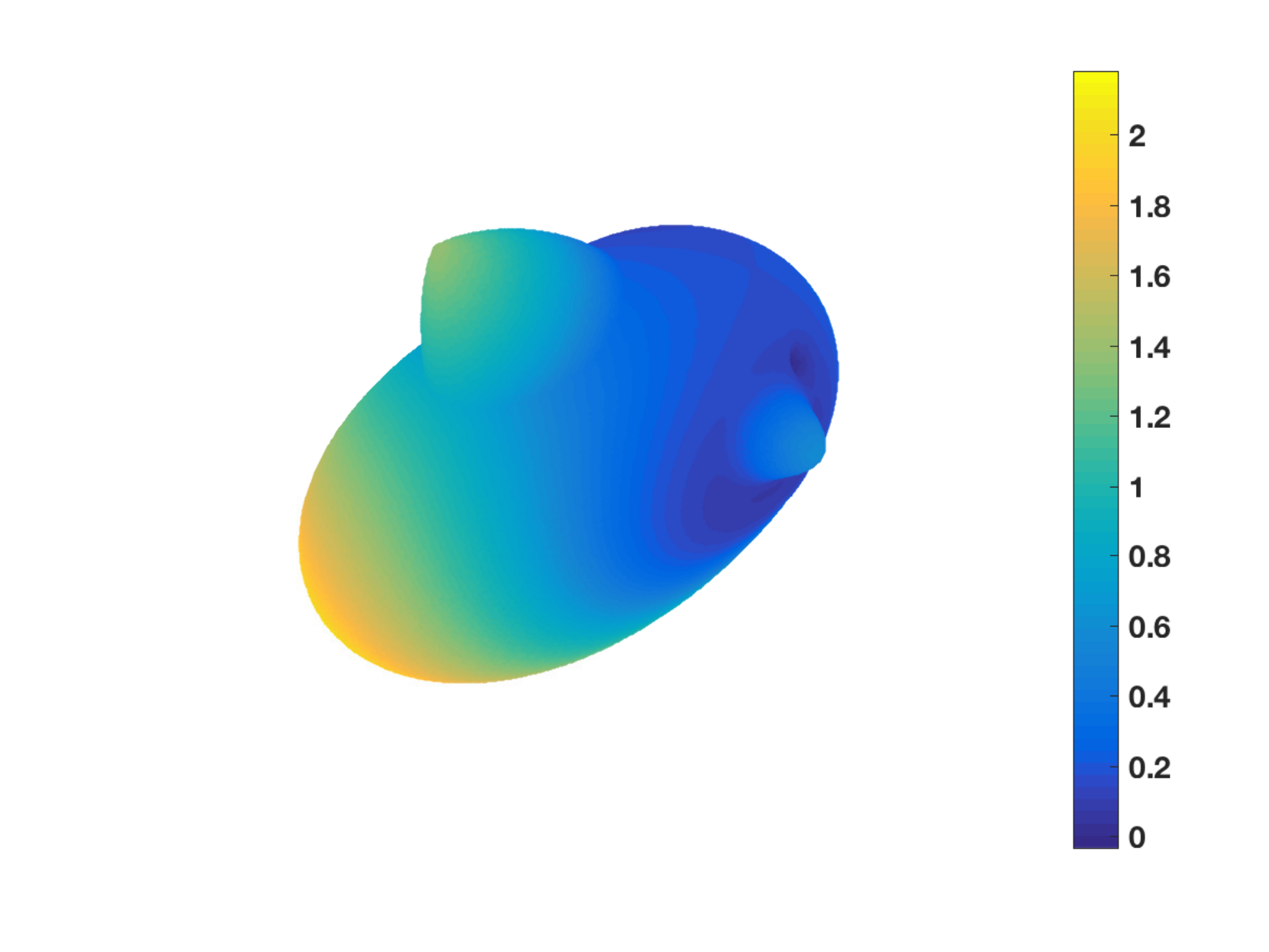}
   }

\subfloat[]{
   \includegraphics[width=0.50\textwidth]{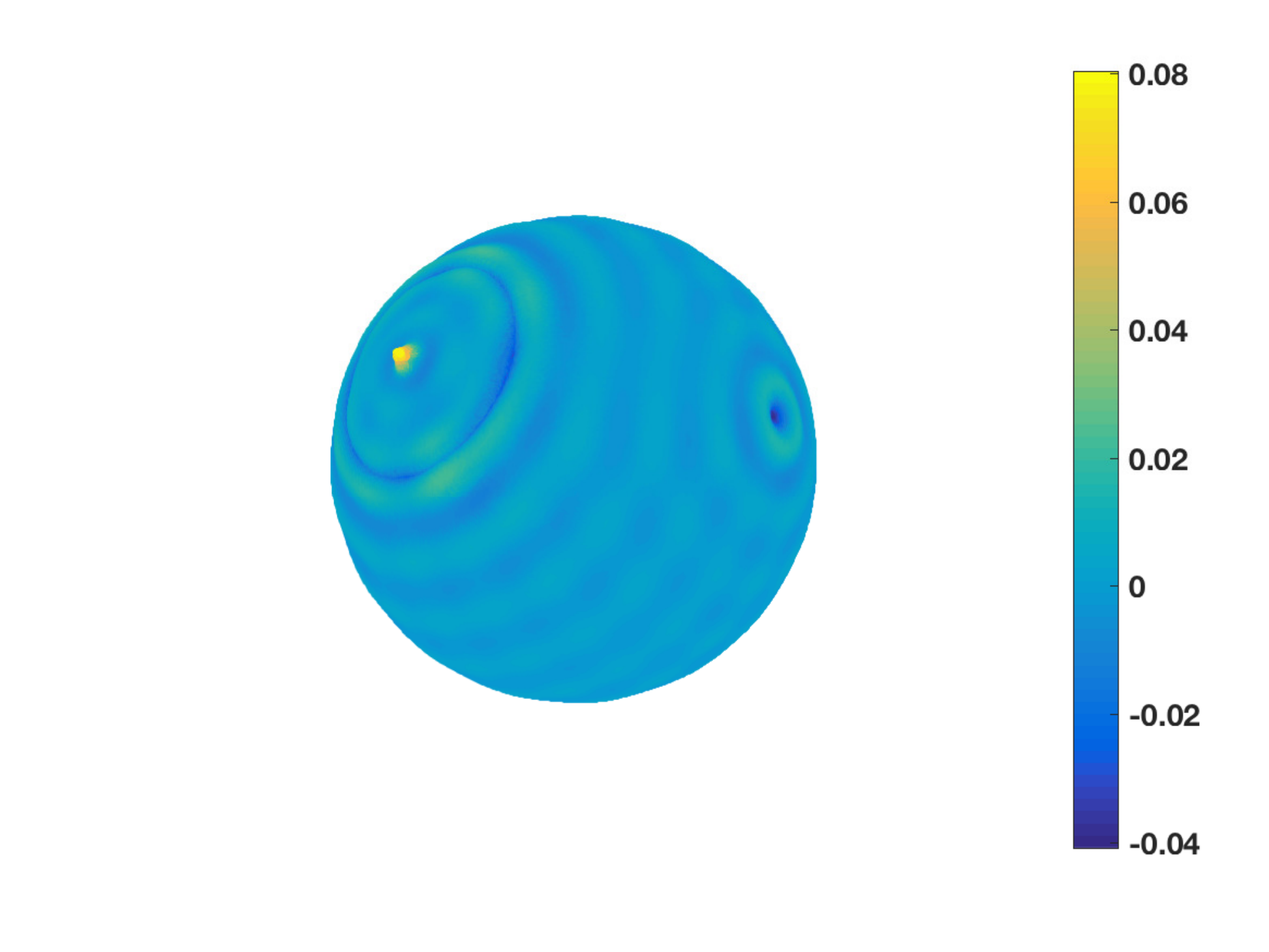}
}
\subfloat[]{
   \includegraphics[width=0.50\textwidth]{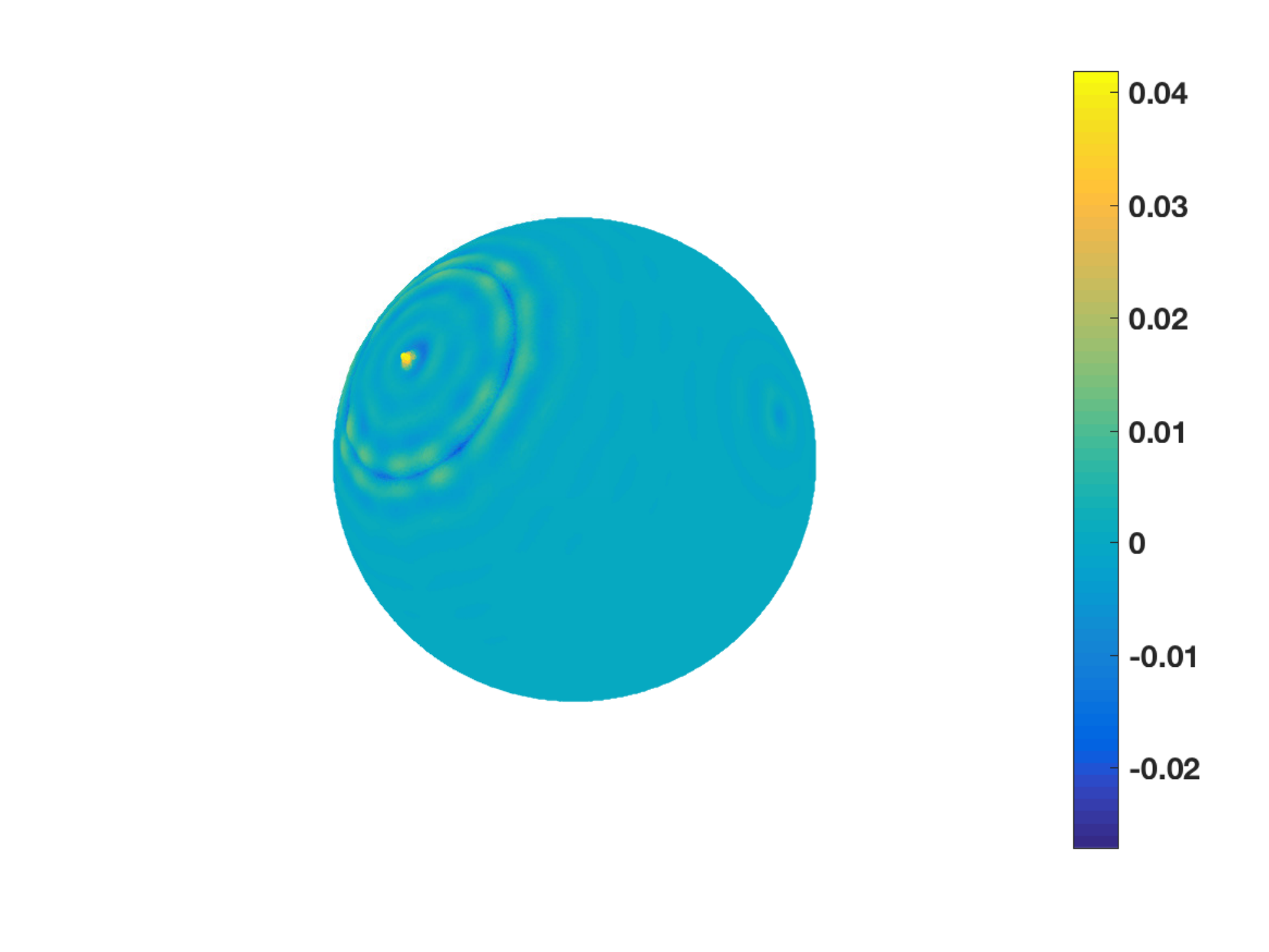}
}
\caption{Top: function $f_3(1,\theta,\varphi)$. Left: hyperinterpolation error. Right: filtered hyperinterpolation error. $K=L=20$ is used.}
\end{figure}

We let $K=20$ and $L=20$ and observe the behaviour of the error. For a comparison, we \red{again} employ interpolation at the Chebyshev zeros in the radial direction, and hyperinterpolation on the sphere in the angular direction. We plot the exact functions at the top, the hyperinterpolation errors at the bottom left, the filtered hyperinterpolation error at the bottom right. Fig.~3--5 show the error on the spherical layers ($r=1,1.0005,1.001$), and Fig.~6 shows the error for the radial line $(-1/2r,-1/2r,1/\sqrt{2}r)$ $(r\in[1,1.001])$. We observe that the filtered hyperinterpolation errors are smaller, and more localised.
\begin{figure}[H]%[htbp]%[H]
\centering
\subfloat[]{
   \includegraphics[width=0.50\textwidth]{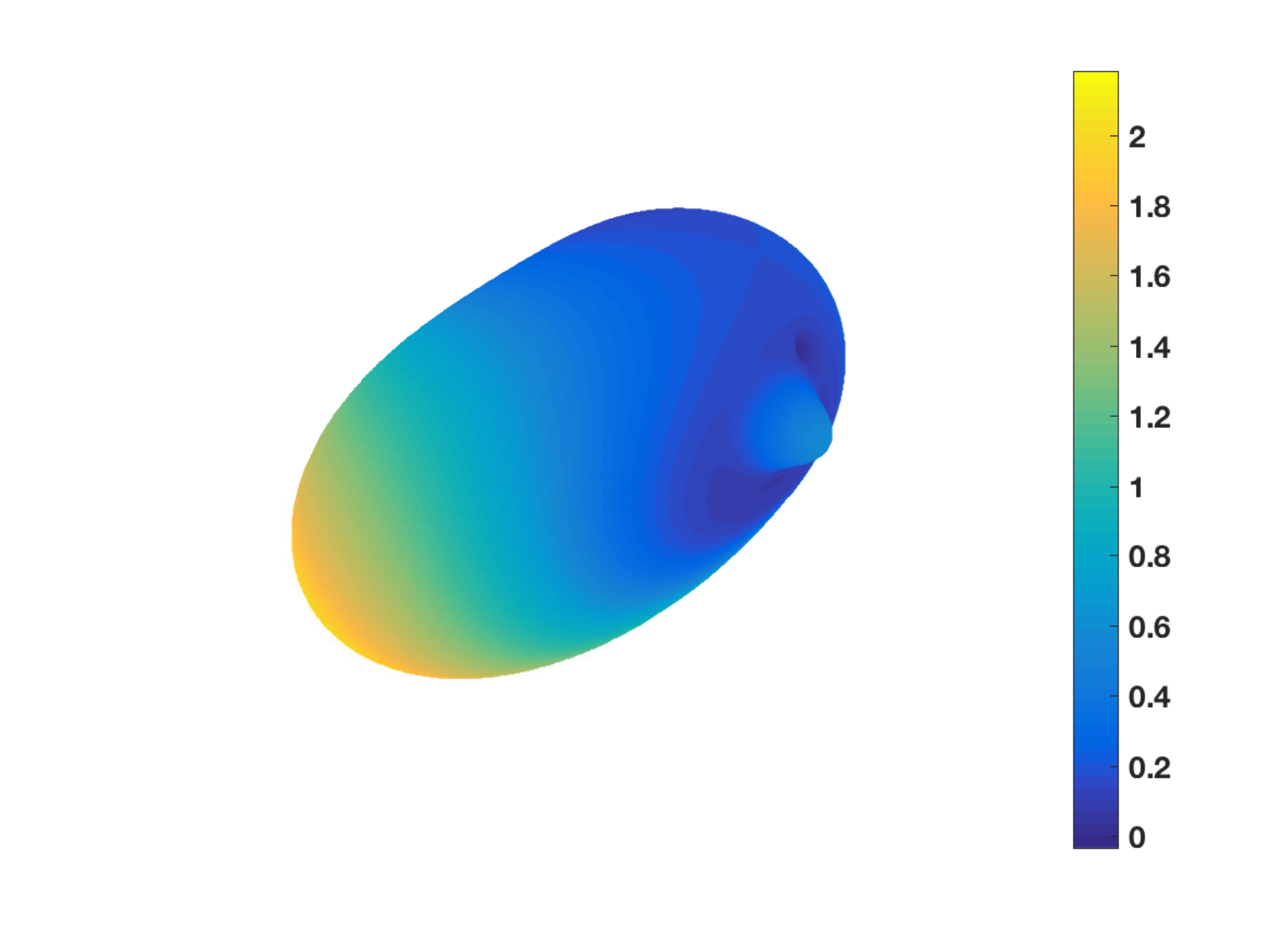}
   }

\subfloat[]{
   \includegraphics[width=0.50\textwidth]{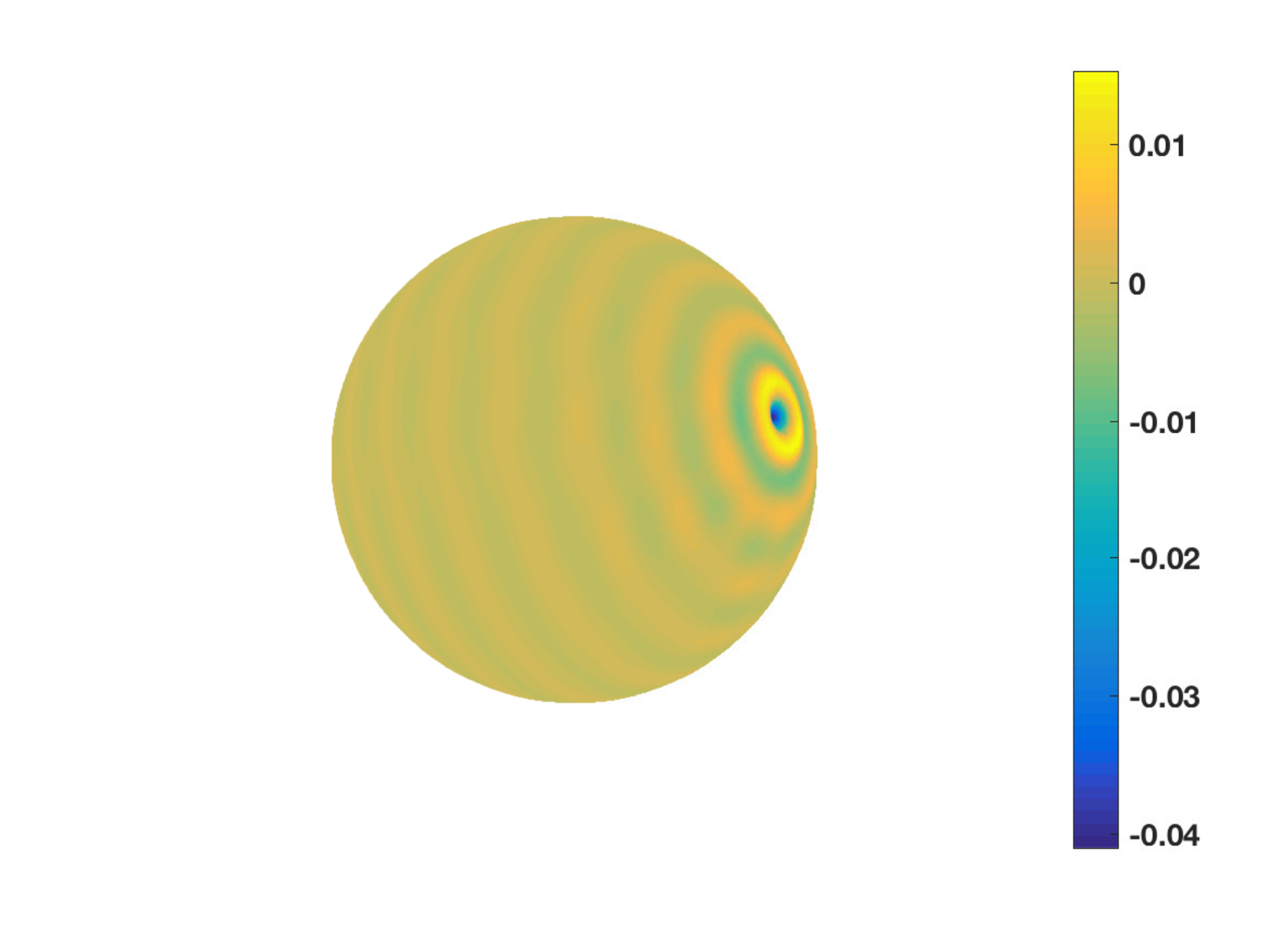}
}
\subfloat[]{
   \includegraphics[width=0.50\textwidth]{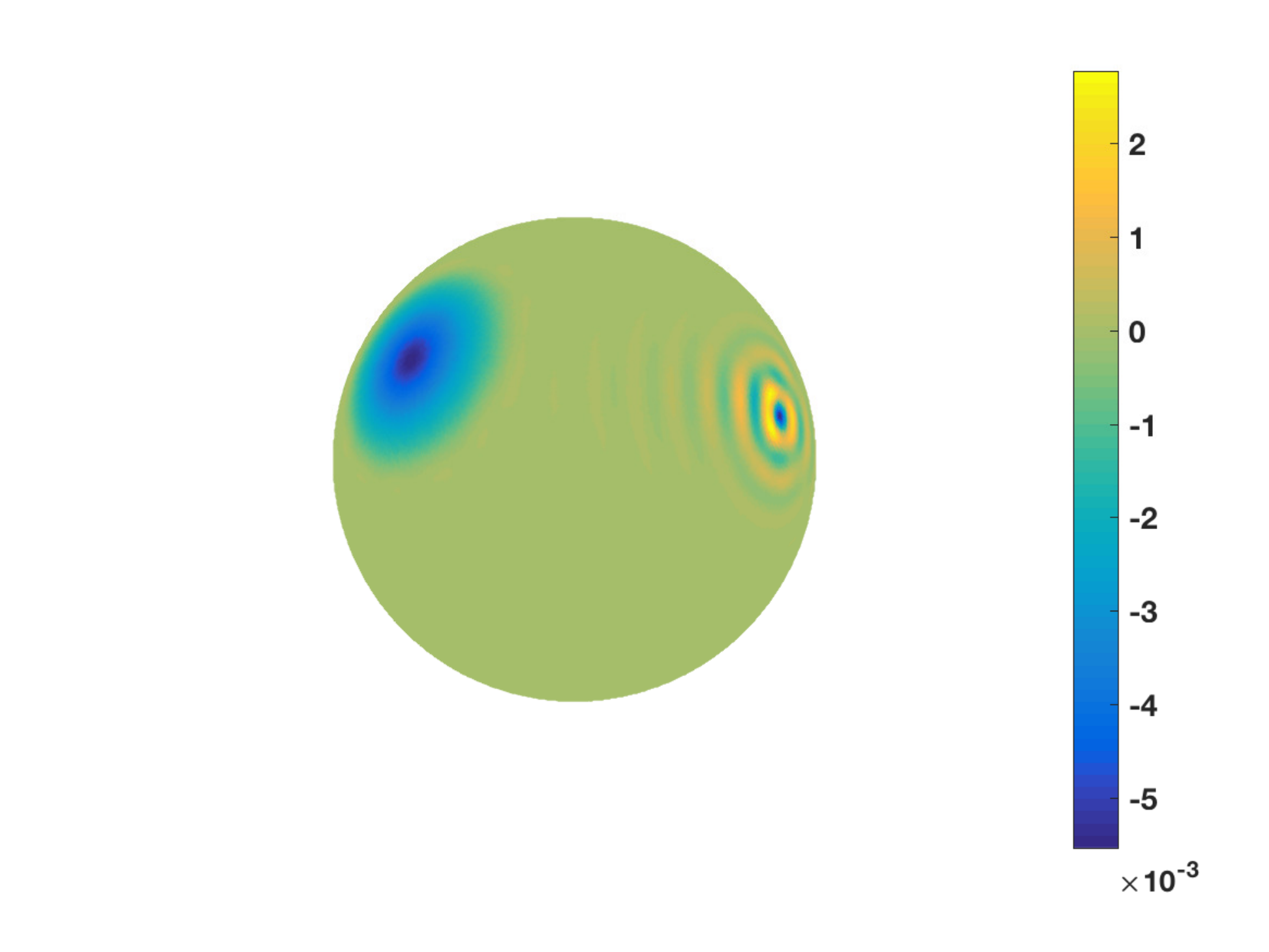}
}
\caption{
Top: function $f_3(1.0005,\theta,\varphi)$. Left: hyperinterpolation error. Right: filtered hyperinterpolation error. $K=L=20$ is used.}
\end{figure}
\begin{figure}[H]%[htbp]%[H]
\centering
\subfloat[]{
   \includegraphics[width=0.50\textwidth]{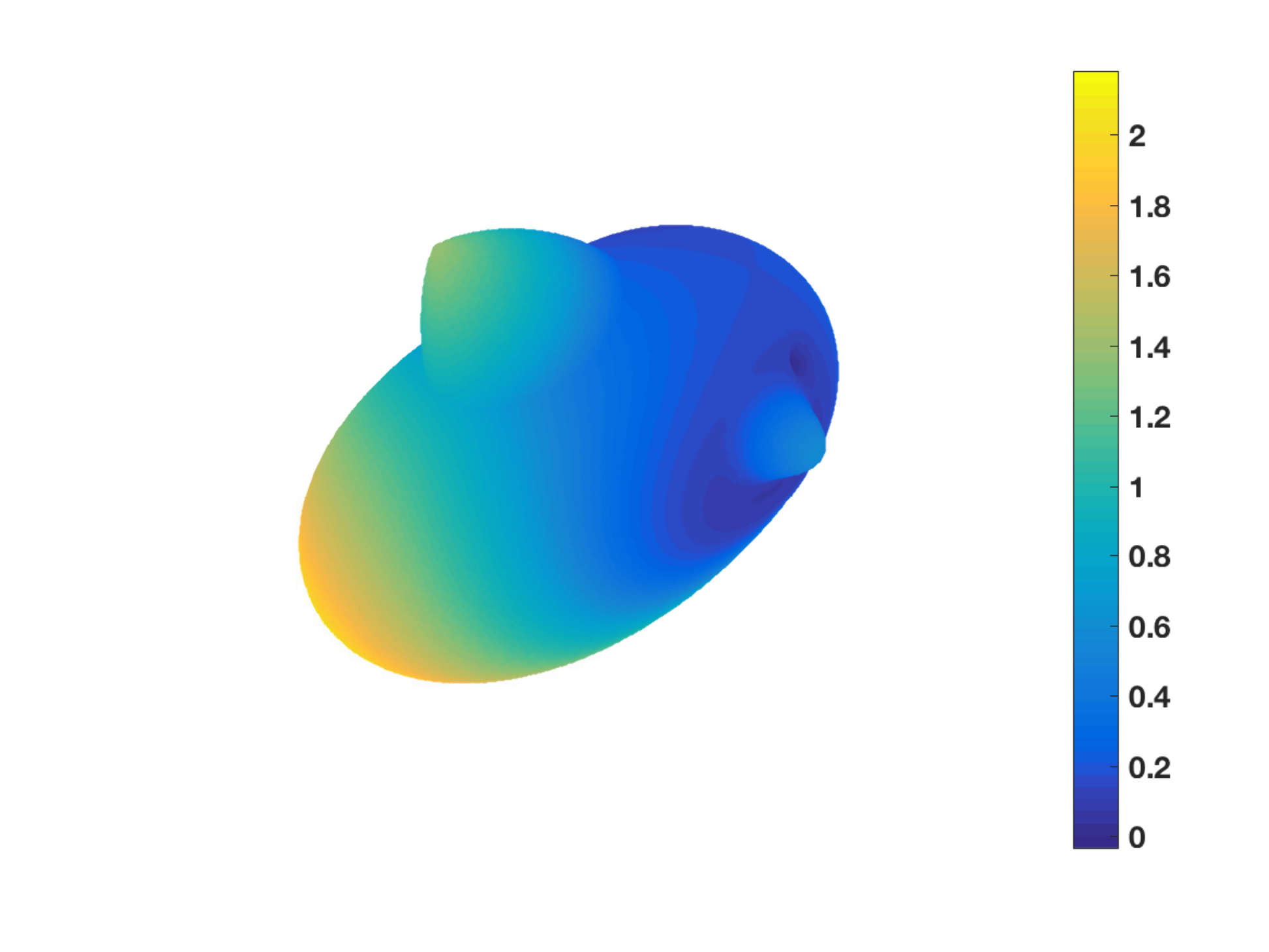}
   }

\subfloat[]{
   \includegraphics[width=0.50\textwidth]{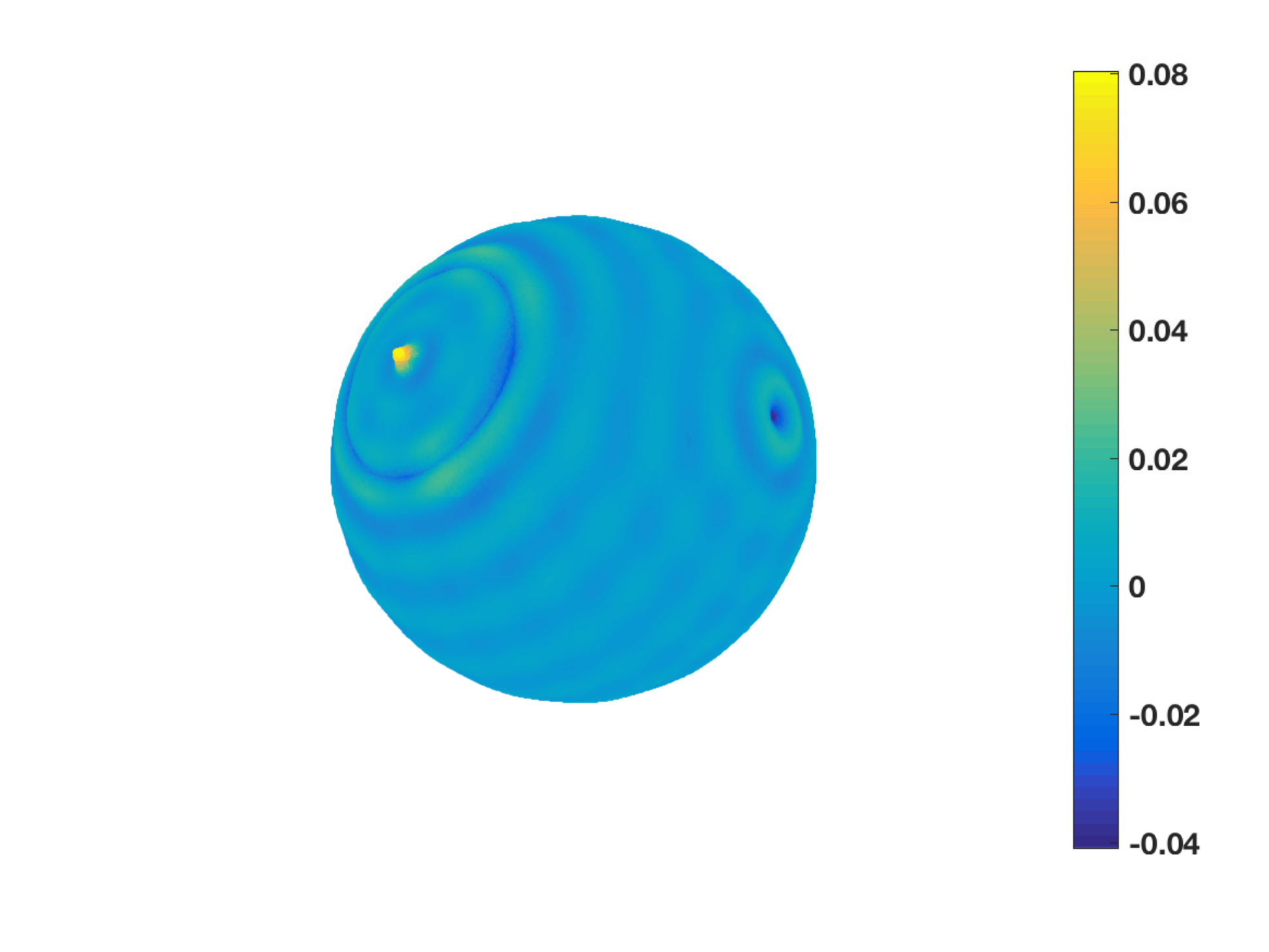}
}
\subfloat[]{
   \includegraphics[width=0.50\textwidth]{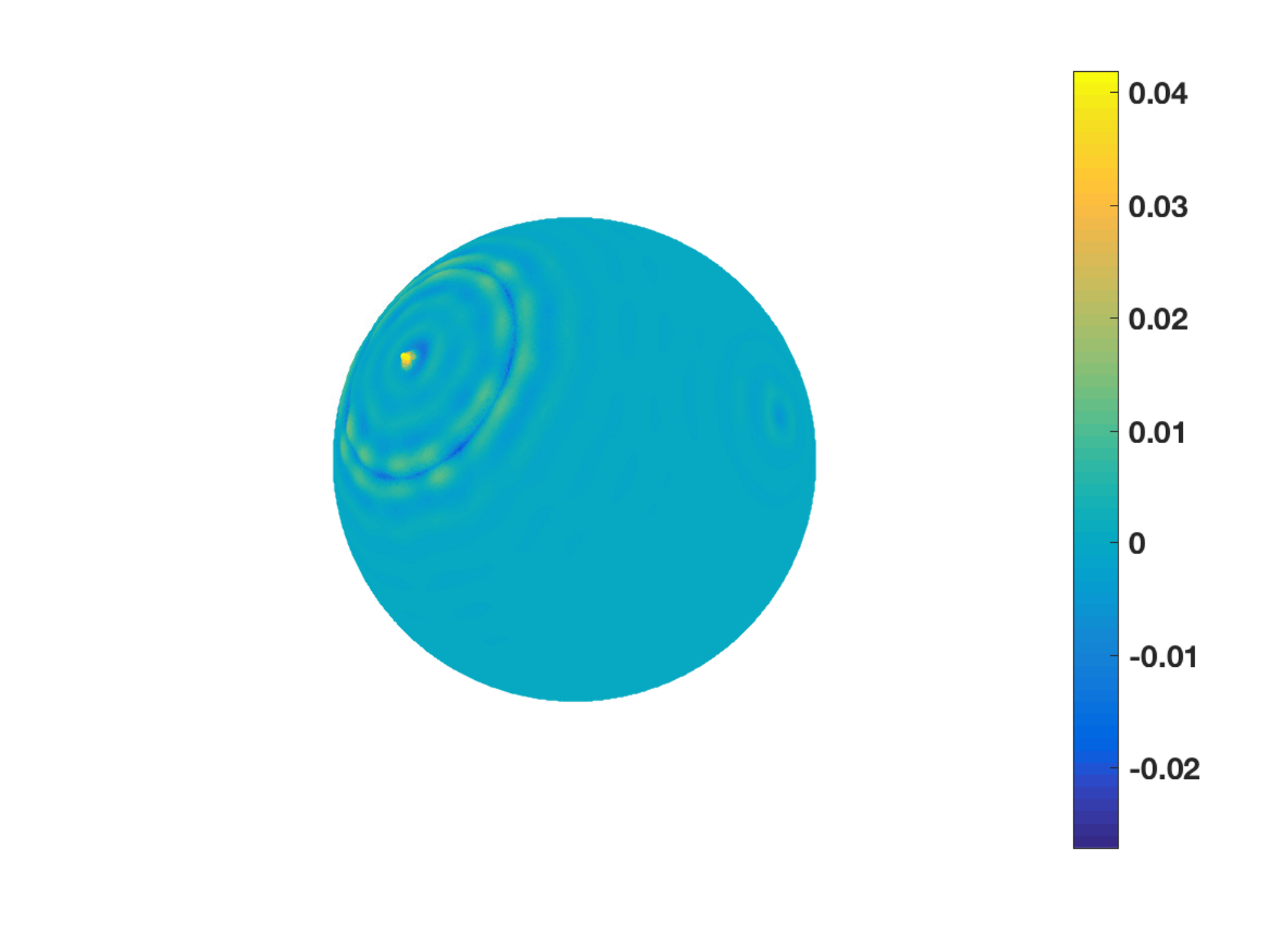}
}
\caption{
Top: function $f_3(1.001,\theta,\varphi)$. Left: hyperinterpolation error. Right: filtered hyperinterpolation error. $K=L=20$ is used.}
\end{figure}
\begin{figure}[H]%[htbp]%[H]
\centering
\subfloat[]{
   \includegraphics[width=0.50\textwidth]{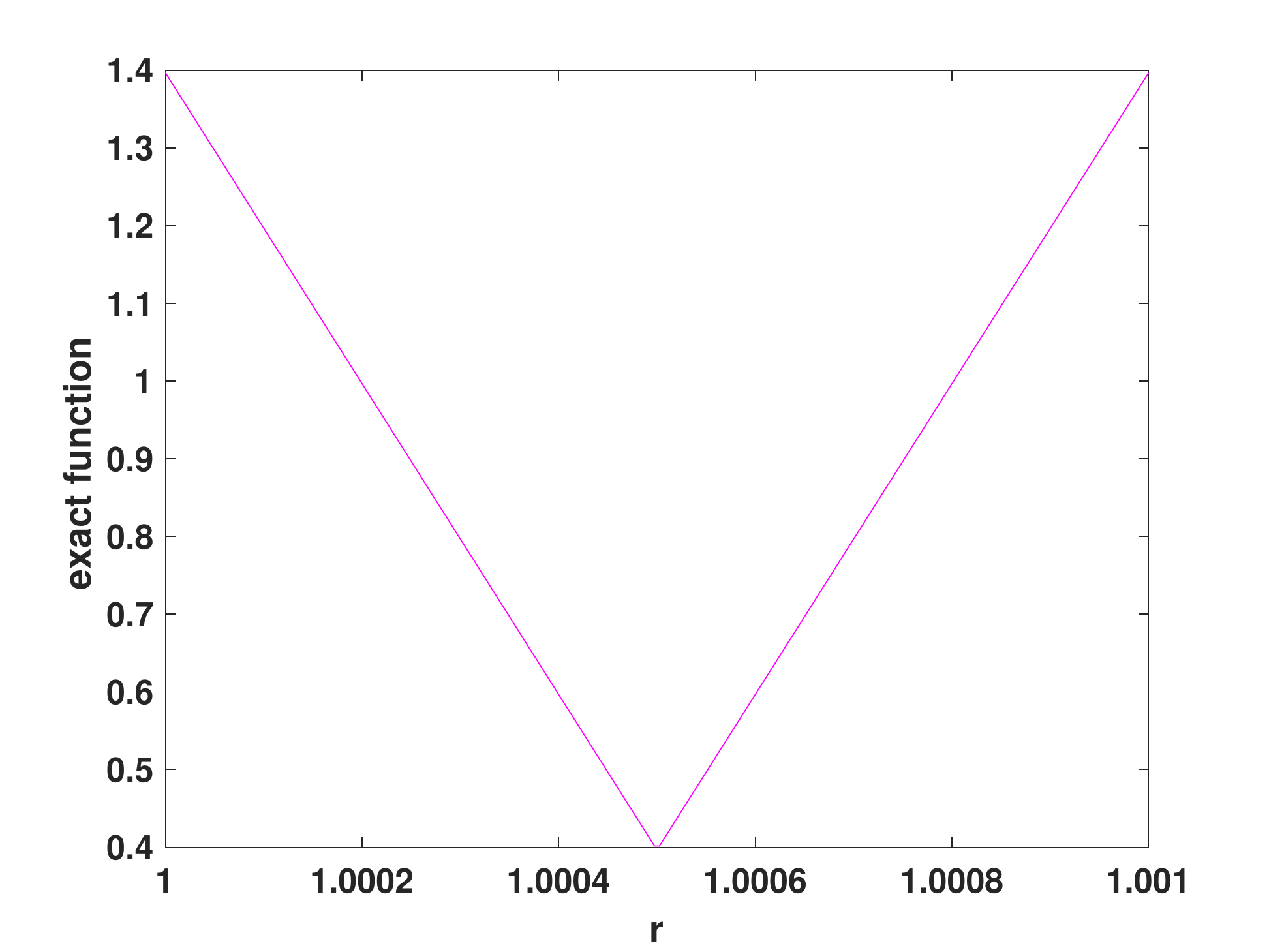}
   }

\subfloat[]{
   \includegraphics[width=0.50\textwidth]{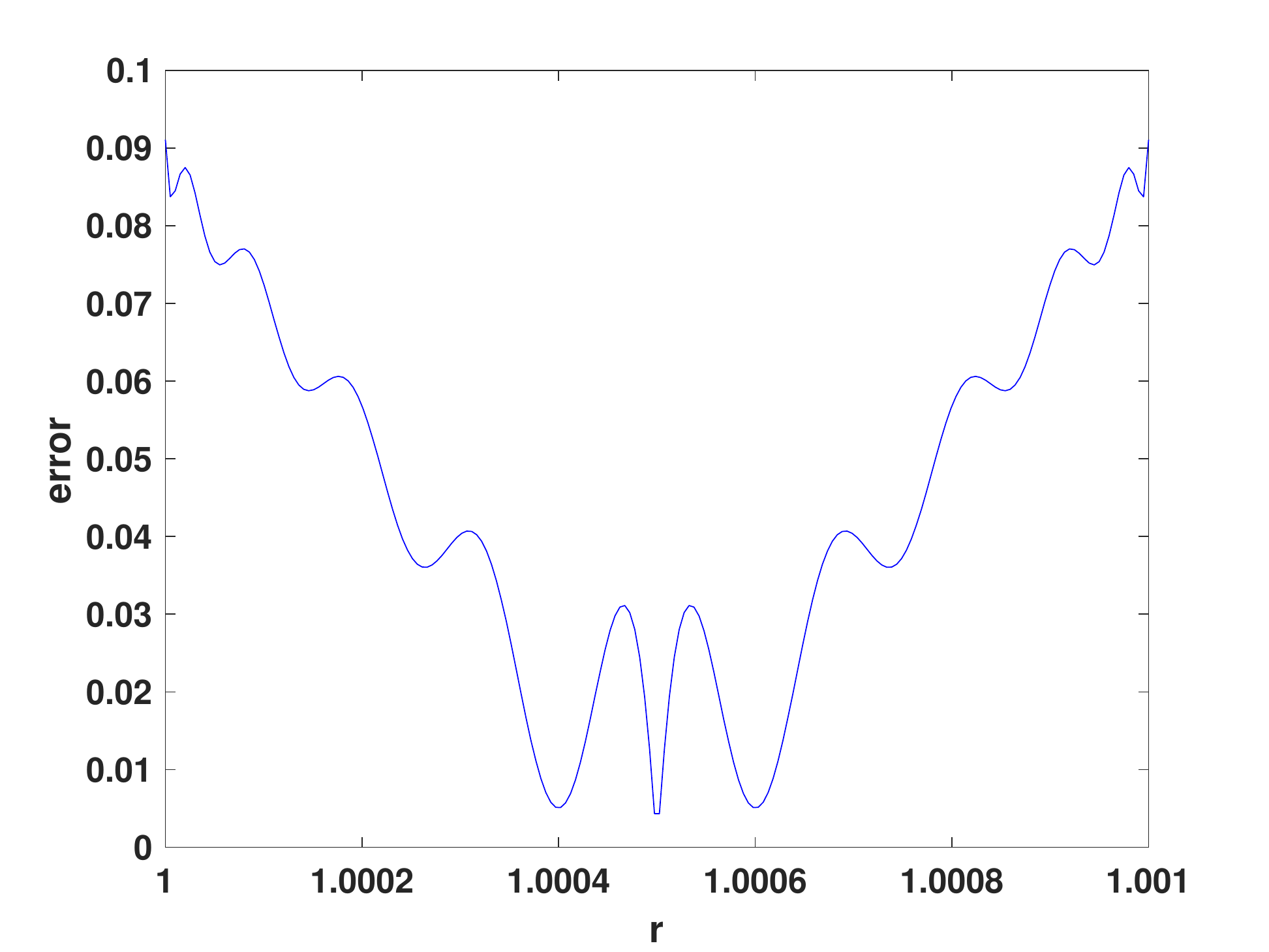}
}
\subfloat[]{
   \includegraphics[width=0.50\textwidth]{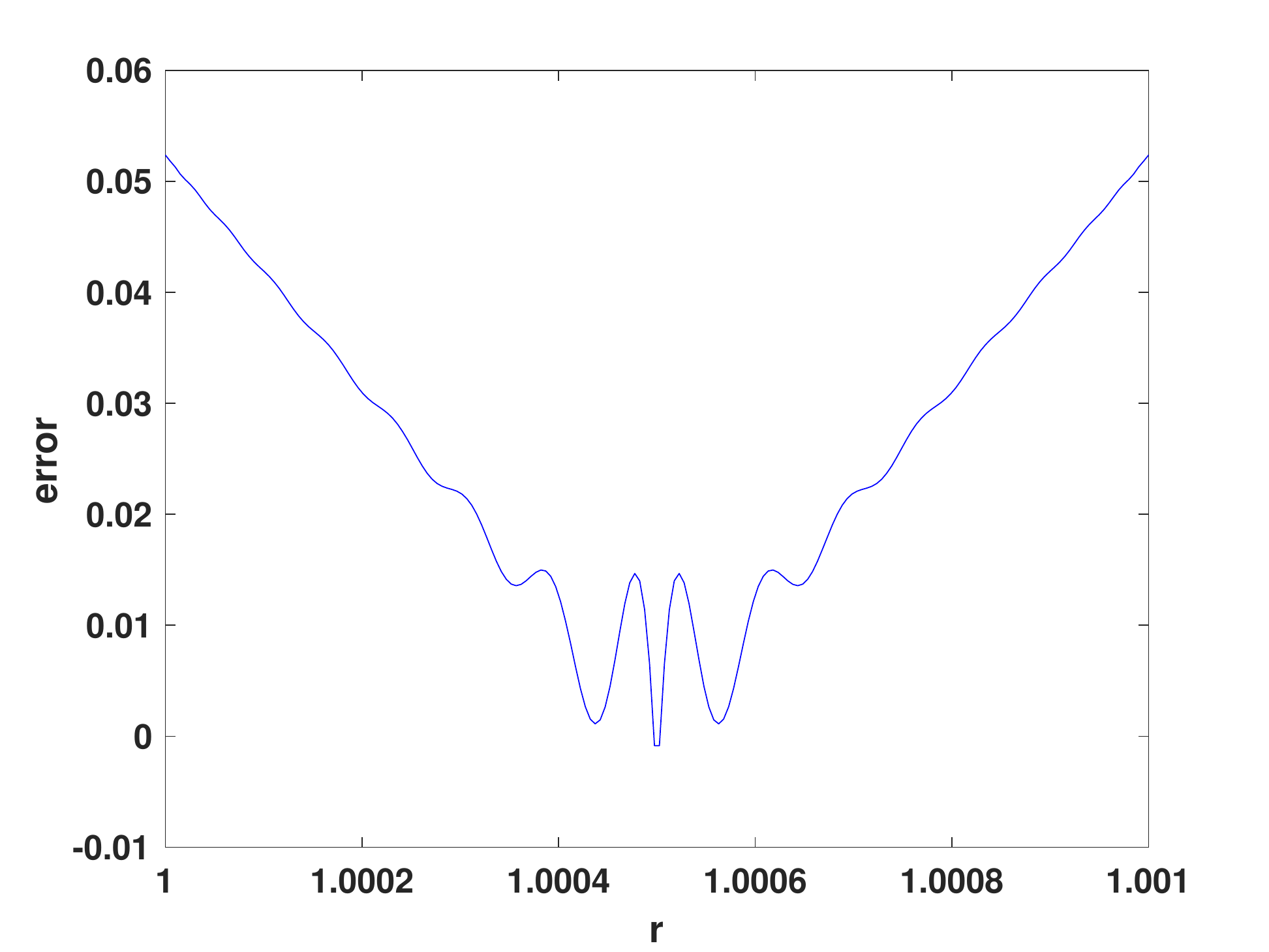}
}
\caption{Top: function $f_3(r,\theta_{\mathrm{c}},\varphi_{\mathrm{c}})$, where $\theta_{\mathrm{c}}$ and $\varphi_{\mathrm{c}}$ are taken so that  $(\sin\theta_{\mathrm{c}}\cos\varphi_{\mathrm{c}},\sin\theta_{\mathrm{c}}\sin\varphi_{\mathrm{c}},\cos\theta_{\mathrm{c}})=(-0.5,-0.5,1/\sqrt{2})$ Left: hyperinterpolation error. Right: filtered hyperinterpolation error.}
\end{figure}
\section{Conclusion}
We proposed a fully discrete filtered polynomial approximation on spherical shells. Our method is based on the filtered hyperinterpolation in the radial and angular directions. We provided an error analysis in terms of the supremum norm by reducing the approximation error to the best polynomial approximations. Numerical results are consistent with the theory.

\section*{Acknowledgements}
I thank Ian H. Sloan and Q. T. Le Gia for helpful discussions. 
\red{I also would like to thank the anonymous reviewer for constructive comments.}
%\bibliography{filtered_abbrev}
\printbibliography
\end{document}

%% file: my_arxiv_header.tex
\usepackage{authblk,setspace}

\usepackage{amsmath} % AMS Math Package 
\usepackage{amsthm} % Theorem Formatting
\usepackage{amssymb}	% Math symbols such as \mathbb
\usepackage{multicol} % Allows for multiple columns
\usepackage{soul} % to use \ul{}, as opposed to \underline{}
\usepackage{cancel}
\usepackage{dsfont,bbm}
\usepackage{empheq,cases}
\usepackage[export]{adjustbox}

\DeclareFontFamily{OT1}{pzc}{}
\DeclareFontShape{OT1}{pzc}{m}{it}{<-> s * [1.10] pzcmi7t}{}
\DeclareMathAlphabet{\mathpzc}{OT1}{pzc}{m}{it}

\usepackage[british]{babel}
\usepackage{csquotes}
\usepackage{calligra}
\DeclareMathAlphabet{\mathcalligra}{T1}{calligra}{m}{n}

\usepackage[low-sup]{subdepth}
\usepackage{hyperref}

\providecommand{\keywords}[1]{\textbf{\textit{Keywords: }} #1}

\let\baraccent=\= % rename builtin command \= to \baraccent
\renewcommand{\=}[1]{\stackrel{#1}{=}} % for putting numbers above =

\newcommand{\norm}[1]{\left\| #1 \right\|}

\newcommand{\ep}{\varepsilon}

\newcommand{\red}[1]{\textcolor{black}{#1}}
\newcommand{\blue}[1]{\textcolor{blue}{#1}}

\newcommand{\bs}[1]{\ensuremath{\boldsymbol{#1}}}

\usepackage{graphicx}
\usepackage{amsmath,amsthm,amssymb,bm,color,comment}
\usepackage{mathrsfs,enumerate}
\usepackage{subfig}
\usepackage{caption}
\usepackage[font=small,format=plain,labelfont=bf,textfont=it]{caption}
\usepackage{url}
\numberwithin{equation}{section}
\definecolor{refkey}{rgb}{0.9451,0.2706,0.4941}
\definecolor{labelkey}{rgb}{0.9451,0.2706,0.4941}
\newtheorem{theorem}{Theorem}[section]
\newtheorem*{theorem*}{Theorem}
\newtheorem{cor}[theorem]{Corollary}

\newtheorem{prop}[theorem]{Proposition}
\theoremstyle{definition}

\newtheoremstyle{myremstyle}
  {\topsep} % Space above
  {1.2\topsep} % Space below
  {} % Body font
  {} % Indent amount
  {\itshape} % Theorem head font {\bfseries}
  {.} % Punctuation after theorem head
  {.5em} % Space after theorem head
  {} % Theorem head spec (can be left empty, meaning `normal')

\theoremstyle{myremstyle} \newtheorem{rem}{Remark}

\usepackage{fancybox,here,mathtools}

\newcommand*{\from}{\colon}

\newcommand{\bbR}{{\mathbb{R}}}
\newcommand{\bbN}{{\mathbb{N}}}

\DeclareMathOperator*{\esssup}{{ess}\sup}

\DeclareMathOperator{\supp}{supp}

\newcommand{\PolyRad}[1]{\mathbb{P}_{#1}(\rinout)}

\newcommand{\dr}{\ensuremath{\,{\rm d}r}}

\newcommand{\dS}{\ensuremath{\,{\rm d}{S}}}

\newcommand{\dmur}{\ensuremath{\,\mathrm{d}\mu} }

\newcommand{\intr}{\ensuremath{\int_{\rin}^{\rout}}}
\newcommand{\shell}{\ensuremath{\mathbb{S}_\varepsilon}}
\newcommand{\interior}[1]{%
  {\mkern0mu#1}^{\mathrm{o}}%
}
\newcommand{\intshell}{\interior{\mathbb{S}}_{\varepsilon}}

\newcommand{\Ltmu}{\ensuremath{L^{2}_{\mu}}}

\newcommand{\Ltmushell}{\ensuremath{{L^{2}_{\mu}(\intshell)}}}

\newcommand{\Qrad}{\ensuremath{{Q_{\mathrm{rad}}}}}

\newcommand{\VKconti}{\ensuremath{\mathcal{R}_{K}}}

\newcommand{\VKL}[1]{\ensuremath{{V}_{KL {#1}}}}

\newcommand{\VK}{\ensuremath{R_{K}}}
\newcommand{\VN}{\ensuremath{{V}^{\mathrm{ang}}_{N}}}
\newcommand{\VL}{\ensuremath{{A}_{L}}}

\newcommand{\sumLwithm}{\ensuremath{\sum_{\ell=0}^{L}\sum_{m=-\ell}^\ell}}

\newcommand{\Ka}{\ensuremath{{\overline{K}(a)}}}
\newcommand{\Lb}{\ensuremath{{\overline{L}(b)}}}

\newcommand{\rin}{\ensuremath{{r_{\mkern0mu\mathrm{in}}}}}
\newcommand{\rout}{\ensuremath{{r_{\mkern-2.5mu\mathrm{out}}}}}
\newcommand{\rinout}{\ensuremath{[\rin,\rout]}}
\newcommand{\weightrad}{\ensuremath{w^{\mathrm{rad}}}}
\newcommand{\weightang}{\ensuremath{w^{\mathrm{ang}}}}

\newcommand{\hrad}{\ensuremath{h^{\mathrm{rad}}\left(\frac{k}K\right)}}
\newcommand{\hang}{\ensuremath{h^{\mathrm{ang}}\left(\frac{\ell}L\right)}}

\newcommand{\innprod}[2]{\left\langle #1, #2 \right\rangle}

\newcommand{\discipS}[2]{\ensuremath{\left\langle #1, #2 \right\rangle^{\mathrm{ang}}_{Q_{\mathrm{ang}}}}}
\newcommand{\discipAng}[3]{\ensuremath{\left\langle #1, #2 \right\rangle^{\mathrm{ang}}_{Q_{\mathrm{ang}}(#3)}}}

\usepackage[
	style=numeric,
	backend=bibtex,
	backref=true,
	sorting=nyt,
	maxbibnames=99,
	firstinits=true,
	doi=false,
	isbn=false%,
	]{biblatex}

\DefineBibliographyStrings{english}{%
  backrefpage = {page},% originally "cited on page"
  backrefpages = {pages},% originally "cited on pages"
}

\renewbibmacro{in:}{}

\DeclareFieldFormat{postnote}{#1}
\DeclareFieldFormat{multipostnote}{#1}
	
\AtEveryBibitem{
\clearfield{month}
\clearfield{number}
\clearlist{location}
}
\AtEveryBibitem{\clearlist{language}} % clears language
\DeclareFieldFormat[article]{title}{#1}
\AtEveryBibitem{%
  \ifentrytype{book}{%
    \clearfield{pages}%
  }{%
  }%
}
\AtEveryBibitem{%
\ifentrytype{book}{
    \clearfield{url}%
}{}
\ifentrytype{article}{
    \clearfield{url}%
}{}
\ifentrytype{incollection}{
    \clearfield{url}%
}{}
}